\documentclass[10pt]{article}

\usepackage{latexsym,amsfonts}
\usepackage{amssymb,amsthm,upref,amscd}
\usepackage[T1]{fontenc}
\usepackage{times}

\newtheorem{Theorem}{Theorem}[section]

\newtheorem{lemma}[Theorem]{Lemma}

\newtheorem{corollary}[Theorem]{Corollary}
\newtheorem{remark}[Theorem]{Remark}
\newtheorem{Open Problem}[Theorem]{Open Problem}

\makeatletter \@addtoreset{equation}{section} \makeatother

\begin{document}

\title{\bf Existence of positive  solutions for a critical  nonlinear
Schr\"odinger equation  with vanishing or coercive
 potentials
 }

\author{Shaowei Chen  \thanks{ Tel.: +86 15059510687. E-mail address:
chensw@amss.ac.cn (S. Chen).
 } \\ \\
\small School of Mathematical Sciences, Huaqiao
University, \\
\small Quanzhou  362021, P.R. China\\
 }

\date{}
\maketitle

\begin{minipage}{13cm}
{\small {\bf Abstract:} In this paper we investigate the existence
of the positive solutions for the following
 nonlinear Schr\"odinger  equation
$$
-\triangle u+V(x)u=K( x)|u|^{p-2}u\ \mbox{in}\ \mathbb{R}^N
$$
where $V(x)\sim a|x|^{-b}$ and $K(x)\sim \mu|x|^{-s}$ as
$|x|\rightarrow\infty$ with $0<a,\mu<+\infty$, $b<2,$ $ b\neq0$,
 $ 0<\frac{s}{b}<1$ and $p=2(N-2s/b)/(N-2).$

}
\end{minipage}

\section{Introduction and statement of results}\label{diyizhang}
In this paper, we consider the following semilinear elliptic
equation
\begin{eqnarray}\label{h777arrrrrr}
-\triangle u+V(x)u=K( x)|u|^{p-2}u\ \mbox{in}\ \mathbb{R}^N.
\end{eqnarray}
where $N\geq 3$. The exponent\begin{eqnarray}\label{lulugenjiuy}
p=2(N-\frac{2s}{b})/(N-2)
\end{eqnarray} with the real numbers   $b$ and $s$
satisfying
\begin{eqnarray}\label{kkh77tygh}
b<2,\ b\neq0,\ 0<\frac{s}{b}<1. \end{eqnarray} By this definition,
$2<p<2^*:=2N/(N-2).$

 With respect to the functions $V$ and $K$, we assume
\begin{description}
\item{$(\bf{A_1}).$}
 $V,K\in C(\mathbb{R}^N)$. For every $
 x\in\mathbb{R}^N$,
 $V(x)>0$ and $K(x)>0$.
  \item{$(\bf{A_2}).$} There exist $0<a<\infty$ and $0<\mu<\infty$ such that
\begin{eqnarray}\label{jjgh7760oi}
\lim_{|x|\rightarrow\infty}|x|^bV(x)=a \ \mbox{and}\
\lim_{|x|\rightarrow\infty}|x|^sK(x)=\mu.
\end{eqnarray}
 \end{description}

A typical example for Eq. (\ref{h777arrrrrr}) with $V$ and $K$
 satisfying $(\bf{A_1})$ and $(\bf{A_2})$ is the equation
 \begin{eqnarray}\label{iifynmcm}
 -\triangle u+\frac{a}{(1+|x|)^b}u=\frac{\mu}{(1+|x|)^s}|u|^{p-2}u\
 \mbox{in}\ \mathbb{R}^N
 \end{eqnarray}
When $0<b<2,$ the potentials are vanishing at infinity and when
$b<0,$ the potentials are coercive.

  Eq.(\ref{h777arrrrrr})
arises in various applications, such as chemotaxis, population
genetics, chemical reactor theory, and the study of standing wave
solutions of certain nonlinear Schr\"odinger equations. Therefore,
they  have received growing attention in recent years (one can
see, e.g., \cite{BGM}, \cite{BL}, \cite{Co}, \cite{P},  \cite{SW}
and \cite{SWW} for reference).

 Under the above  assumptions, Eq.(\ref{h777arrrrrr}) has a  natural variational structure.
 For an open subset $\Omega$ in $\mathbb{R}^N,$ let $C^\infty_0(\Omega)$ be  the collection of
smooth functions with compact support set in $\Omega$. Let $E$ be
the completion of $C^\infty_0(\mathbb{R}^N)$ with respect to the
inner product
$$(u,v)_E=\int_{\mathbb{R}^N}\nabla u\nabla
vdx+\int_{\mathbb{R}^N}V(x)uvdx.$$ From the assumptions $\bf
(A_1)$ and $\bf (A_2)$, we deduce that
$$(\int_{\mathbb{R}^N}\frac{|u|^{2}}{(1+|x|)^b}dx)^{1/2}\ \mbox{and}\ (\int_{\mathbb{R}^N}V(x)|u|^{2}dx)^{1/2}$$ are two equivalent
norms in the space
$$L^2_V(\mathbb{R}^N)=\{u \ \mbox{is measurable in}\ \mathbb{R}^N
\ |\ \int_{\mathbb{R}^N}V(x)|u|^2dx<+\infty\}.$$ Therefore, there
exists $B_1>0$ such that
$$(\int_{\mathbb{R}^N}\frac{|u|^{2}}{(1+|x|)^b}dx)^{1/2}\leq B_1(\int_{\mathbb{R}^N}V(x)|u|^{2}dx)^{1/2}.$$
Moreover, the assumptions $\bf (A_1)$ and $\bf (A_2)$ imply that
there exists $B_2>0$ such that
$$K(x)\leq B_2(1+|x|)^{-s},\ \forall x\in\mathbb{R}^N.$$  Then by the H\"older
and  the Sobolev inequalities (see, e.g., \cite[Theorem
1.8]{Willem}), we have, for every $u\in C^\infty_0(\mathbb{R}^N)$,
 \begin{eqnarray}\label{hhfgtr5ftz}
(\int_{\mathbb{R}^N} K(x)|u|^pdx)^{\frac{1}{p}}
&\leq&C(\int_{\mathbb{R}^N}\frac{|u|^p}{(1+|x|)^s}dx)^{\frac{1}{p}}\nonumber\\
&=&C(\int_{\mathbb{R}^N}\frac{|u|^{\frac{2s}{b}}}{(1+|x|)^s}\cdot|u|^{p-\frac{2s}{b}}dx)^{\frac{1}{p}}\nonumber\\
&\leq&C(\int_{\mathbb{R}^N}\frac{|u|^{2}}{(1+|x|)^b}dx)^{\frac{s}{pb}}(\int_{\mathbb{R}^N}|u|^{2^*}dx)^{\frac{1}{p}(1-\frac{s}{b})}\nonumber\\
&\leq&C(\int_{\mathbb{R}^N}\frac{|u|^{2}}{(1+|x|)^b}dx)^{\frac{s}{pb}}(\int_{\mathbb{R}^N}|\nabla
u|^{2}dx)^{\frac{2^*}{2p}(1-\frac{s}{b})}\nonumber\\
&=&C(\int_{\mathbb{R}^N}\frac{|u|^{2}}{(1+|x|)^b}dx)^{\frac{1}{2}\cdot\frac{2s}{pb}}(\int_{\mathbb{R}^N}|\nabla
u|^{2}dx)^{\frac{1}{2}\cdot(1-\frac{2s}{pb})}\nonumber\\
&\leq&C(\int_{\mathbb{R}^N}V(x)|u|^{2}dx)^{\frac{1}{2}\cdot\frac{2s}{pb}}(\int_{\mathbb{R}^N}|\nabla
u|^{2}dx)^{\frac{1}{2}\cdot(1-\frac{2s}{pb})},\nonumber
 \end{eqnarray}
where $C>0$ is a constant independent of $u$. It follows that
there exists a constant $C'>0$ such that
 $$(\int_{\mathbb{R}^N}K(x)|u|^pdx)^{1/p}\leq C'(\int_{\mathbb{R}^N}|\nabla u|^2dx)^{1/2}
 +C'(\int_{\mathbb{R}^N}V(x)|u|^2dx)^{1/2}.$$
This  implies that  $E$ can be embedded continuously into the
weighted $L^p-$space
$$L^p_K(\mathbb{R}^N)=\{u \ \mbox{is measurable in}\ \mathbb{R}^N
\ |\ \int_{\mathbb{R}^N}K(x)|u|^pdx<+\infty\}.$$ Then the
functional
$$\Phi(u)=\frac{1}{2}||u||^2_E-\frac{1}{p}\int_{\mathbb{R}^N}K(x)|u|^pdx,\ u\in E$$
is well defined in $E$. And it is easy to check that $\Phi$ is
 a $C^2$ functional and  the critical points of $\Phi$ are solutions of (\ref{h777arrrrrr}) in $E$.

  In a recent paper \cite{AS}, Alves
and Souto proved that the space $E$ can be embedded compactly into
$L^p_K(\mathbb{R}^N)$ if $0<b<2$ and $2(N-2s/b)/(N-2)<p<2^*$ and
$\Phi$ satisfies Palais-Smale condition consequently. Then by
using the mountain pass theorem, they obtained a nontrivial
solution for Eq.(\ref{h777arrrrrr}). Unfortunately, when
$p=2(N-2s/b)/(N-2)$, the embedding of $E$ into
$L^p_K(\mathbb{R}^N)$ is not compact and
 $\Phi$  satisfies no longer Palais-Smale condition.
 Therefore, the "standard" variational methods fail in this case.
  From this point of view, $p=2(N-2s/b)/(N-2)$ should
be seen as a kind of critical exponent for Eq.(\ref{h777arrrrrr}).
If the potentials $V$ and $K$ are restricted to the class of
radially symmetric functions, "compactness" of such a kind is
regained and "standard" variational approaches work (see \cite{SW}
and \cite{SWW}). But this method does not seem to apply to the
more general equation (\ref{h777arrrrrr}) where $K$ and $V$ are
non-radially symmetric functions.

It is not easy to deal with Eq. (\ref{h777arrrrrr}) directly
because there are no known approaches can be used directly to
overcome the difficulty brought by the loss of compactness.
However, in this paper, through an interesting transformation, we
find an equivalent equation for Eq. (\ref{h777arrrrrr}) (see Eq.
(\ref{o8867yyttggd}) in Section \ref{uuat66atrrsf}). This equation
has the advantages that its Palais-Smale sequence can be
characterized precisely through the concentration-compactness
principle (see Theorem \ref{uut6ghvbbo}) and it possesses partial
compactness (see Corollary \ref{kkg77t6yy}). By means of these
advantages, a positive solution for this equivalent equation and
then a corresponding positive solution for Eq. (\ref{h777arrrrrr})
are obtained.

 Before
to state our main result, we need to give some definitions.

 Let \begin{eqnarray}\label{runacdrf65}
V_*(x)=|x|^{\frac{2b}{2-b}}V(|x|^{\frac{b}{2-b}}x)+C_b|x|^{-2},\end{eqnarray}
where \begin{eqnarray}\label{jjgf7745ttr}
C_b=\frac{b}{4}(1-\frac{b}{4})(N-2)^2\end{eqnarray} and
\begin{eqnarray}\label{2runacdrf65}
K_*(x)=|x|^{\frac{2s}{2-b}}K(|x|^{\frac{b}{2-b}}x).
\end{eqnarray}
 Let
$H^1(\mathbb{R}^N)$ be the the Sobolev space endowed  with the
norm and the inner product
$$||u||=(\int_{\mathbb{R}^N}|\nabla u|^2dx+\int_{\mathbb{R}^N}u^2dx)^{1/2} \ \mbox{and}\
( u,v) =\int_{\mathbb{R}^N}(\nabla u\cdot\nabla v+uv)dx$$
respectively and $L^p(\mathbb{R}^N)$ be the function space
consisting of the  functions  on $\mathbb{R}^N$ that are
$p-$integrable. Since $2<p<2^*,$ $H^1(\mathbb{R}^N)$ can be
embedded continuously into   $L^p(\mathbb{R}^N).$
 Therefore, the infimum
\begin{eqnarray}\label{olfgcbb5er}
\inf_{v\in
H^1(\mathbb{R}^N)\setminus\{0\}}\frac{\int_{\mathbb{R}^N}|\nabla
v|^2dx+a\int_{\mathbb{R}^N}v^2dx}{(\int_{\mathbb{R}^N}|v|^pdx)^{2/p}}>0.
\end{eqnarray}
We denote this infimum by $S_p.$

Our main result   reads as follows:

 \begin{Theorem}\label{ll7trggfgf}
Under the assumptions  $\bf (A_1)$ and $\bf (A_2)$, if $b,s$ and
$p$ satisfy (\ref{kkh77tygh}) and (\ref{lulugenjiuy}) and
\begin{eqnarray}\label{ofbbv66ft}
&&\inf_{u\in
H^1(\mathbb{R}^N)\setminus\{0\}}\frac{\int_{\mathbb{R}^N}|\nabla
u|^2dx+(\frac{b^2}{4}-b)\int_{\mathbb{R}^N}\frac{|x\cdot\nabla
u|^2}{|x|^2}dx+\int_{\mathbb{R}^N}V_*(x)|u|^2dx}{(\int_{\mathbb{R}^N}K_*(x)|u|^pdx)^{2/p}}\nonumber\\
&<&(1-b/2)^{\frac{p-2}{p}}\mu^{-\frac{2}{p}}S_p,
\end{eqnarray} then Eq. (\ref{h777arrrrrr}) has a positive solution $u\in E$.
 \end{Theorem}
 \begin{remark} We should emphasize
that the condition (\ref{ofbbv66ft}) can be satisfied in many
situations. For $r>0,$ let $R_r=\{x\in\mathbb{R}^N\ |\
r/2<|x|<r\}$ and $H^1_0(R_r)$ be the closure of $C^\infty_0(R_r)$
in $H^1(\mathbb{R}^N)$. Under the assumptions $\bf (A_1)$ and $\bf
(A_2)$, we have
$$\inf_{u\in
H^1_0(R_r)\setminus\{0\}}\frac{\int_{R_r}|\nabla
u|^2dx}{(\int_{R_r}K_*(x)|u|^pdx)^{2/p}}\rightarrow 0,\ \mbox{as}\
r\rightarrow +\infty.$$ Then for any $\epsilon>0$, there exist
$r_\epsilon>0$ and $u_\epsilon\in H^1_0(R_r)\setminus\{0\}$ such
that
$$\frac{\int_{R_r}|\nabla
u_\epsilon|^2dx}{(\int_{R_r}K_*(x)|u_\epsilon|^pdx)^{2/p}}<\epsilon.$$
It follows from this inequality and $\int_{R_r}\frac{|x\cdot\nabla
u_\epsilon|^2}{|x|^2}dx\leq \int_{R_r}|\nabla u_\epsilon|^2dx$
that if $\sup_{R_r}V_*$ is  small enough such
 that
$$\frac{\int_{R_r}V_*(x)|u_\epsilon|^2dx}{(\int_{R_r}K_*(x)|u_\epsilon|^pdx)^{2/p}}<\epsilon,$$
then
\begin{eqnarray}\label{kkf66ftrtqas}
&&\frac{\int_{R_r}|\nabla
u_\epsilon|^2dx+(\frac{b^2}{4}-b)\int_{R_r}\frac{|x\cdot\nabla
u_\epsilon|^2}{|x|^2}dx+\int_{R_r}V_*(x)|u_\epsilon|^2dx}{(\int_{R_r}K_*(x)|u_\epsilon|^pdx)^{2/p}}\nonumber\\
&<&(2+|\frac{b^2}{4}-b|)\epsilon\nonumber
\end{eqnarray}
This implies that (\ref{ofbbv66ft}) is satisfied if $\epsilon$ is
chosen such that
$(2+|\frac{b^2}{4}-b|)\epsilon<(1-b/2)^{\frac{p-2}{p}}\mu^{-\frac{2}{p}}S_p.$
 \end{remark}
\noindent{\bf Notations}: Let $X$ be a Banach Space and
$\varphi\in C^1(X,\mathbb{R}).$ We denote the Fr\'echet derivative
of $\varphi$ at $u$ by $\varphi'(u)$.  The Gateaux derivative of
$\varphi$ is denoted by $\langle \varphi'(u), v\rangle,$ $\forall
u,v\in X.$ By $\rightarrow$ we denote the strong and by
$\rightharpoonup$ the weak convergence. For a function $u$,
$u^{+}$ denotes the functions $\max\{u(x),0\}$. The symbol
$\delta_{ij}$ denotes the Kronecker symbol:
 $\delta_{ij}=\left\{
\begin{array}{l}
1, \ i=j\\
0, \ i\neq j.\\
\end{array} \right.$ We  use $o(h)$ to mean
 $o(h)/|h|\rightarrow 0$ as $|h|\rightarrow
0$.

\section{An equivalent equation for Eq. (\ref{h777arrrrrr})}\label{uuat66atrrsf}
For $x\in\mathbb{R}^N,$ let $y=|x|^{-b/2}x$. To $u,$ a $C^2$
function in $\mathbb{R}^N,$ we associate a function $v$, a $C^2$
function
 in $\mathbb{R}^N\setminus\{0\}$ by the transformation
\begin{eqnarray}\label{lululuggftrr}
u(x)=|x|^{-\frac{b}{4}(N-2)}v(|x|^{-\frac{b}{2}}x_1,\cdots,
|x|^{-\frac{b}{2}}x_N).
\end{eqnarray}
\begin{lemma}\label{jjfgr665t2w}
Under the above assumptions,
\begin{eqnarray}\label{kkg7754r55}
\triangle_x
u(x)=|y|^{-\frac{b(N+2)}{2(2-b)}}\Big(\sum_{i,j=1}^N\frac{\partial}{\partial
y_j}\Big(A_{ij}(y)\frac{\partial v}{\partial
y_i}\Big)-\frac{C_b}{|y|^2}v\Big).
\end{eqnarray}
where
\begin{eqnarray}\label{oo77565ttr}
A_{ij}(y)=\delta_{ij}+(\frac{b^2}{4}-b)\frac{y_iy_j}{|y|^2},\
i,j=1,\cdots,N.
\end{eqnarray}
\end{lemma}
\noindent{\bf Proof.} Let $r=|x|$. By   direct computations,
\begin{eqnarray}\label{kkg77tyg6}
\frac{\partial u}{\partial x_i}
&=&r^{-\frac{b(N-2)}{4}-\frac{b}{2}}\frac{\partial v}{\partial
y_i}-\frac{b}{2}r^{-\frac{b(N-2)}{4}-\frac{b}{2}-2}x_i\sum_{j=1}^N
x_{j}\frac{\partial v}{\partial
y_j}\nonumber\\
&&-\frac{b}{4}(N-2)r^{-\frac{b(N-2)}{4}-2}x_iv
\end{eqnarray}
and
\begin{eqnarray}\label{igu77g6ty}
\frac{\partial^2 u}{\partial x^2_i}
&=&-\frac{bN}{2}r^{-\frac{b(N-2)}{4}-\frac{b}{2}-2}x_i\frac{\partial
v}{\partial y_i}+r^{-\frac{b(N-2)}{4}-b}\frac{\partial^2
v}{\partial
y^2_i}-br^{-\frac{b(N-2)}{4}-b-2}\sum^{N}_{j=1}x_jx_i\frac{\partial^2
v}{\partial y_j\partial
y_i}\nonumber\\
&&+\Big(\frac{b^2}{4}(N-1)+b\Big)r^{-\frac{b(N-2)}{4}-\frac{b}{2}-4}x^2_i\sum^{N}_{j=1}x_j\frac{\partial
v}{\partial y_j}\nonumber\\
&&-\frac{b}{2}r^{-\frac{b(N-2)}{4}-\frac{b}{2}-2}\sum^{N}_{j=1}x_j\frac{\partial
v}{\partial
y_j}+\frac{b^2}{4}r^{-\frac{b(N-2)}{4}-b-4}x^2_i\sum^{N}_{j,k=1}x_jx_k\frac{\partial^2
v}{\partial y_j\partial y_k}\nonumber\\
&&+\frac{b}{4}(N-2)(\frac{b}{4}(N-2)+2)r^{-\frac{b}{4}(N-2)-4}x^2_iv-\frac{b}{4}(N-2)r^{-\frac{b}{4}(N-2)-2}v.\nonumber
\end{eqnarray}
Then
\begin{eqnarray}\label{skkg77tyg6}
\triangle_x u&=&\sum^{N}_{i=1}\frac{\partial^2 u}{\partial
x_i^2}\nonumber\\
&=&r^{-\frac{b(N-2)}{4}-b}\Big\{\triangle_yv+(\frac{b^2}{4}-b)r^{-2}\sum^{N}_{i,j=1}x_ix_j\frac{\partial^2
v}{\partial y_i\partial y_j}\nonumber\\
&&+(\frac{b^2}{4}-b)(N-1)r^{\frac{b}{2}-2}\sum_{i=1}^N
x_i\frac{\partial v}{\partial
y_i}-\frac{b}{4}(1-\frac{b}{4})(N-2)^2r^{b-2}v\Big\}.
\end{eqnarray}
Since $y=|x|^{-b/2}x,$ we have   $r=|y|^{\frac{2}{2-b}}$ and
$x_i=|y|^{\frac{b}{2-b}}y_i,$ $1\leq i\leq N$. Then
\begin{eqnarray}\label{kkg88765te}
&&r^{-2}\sum^{N}_{i,j=1}x_ix_j\frac{\partial^2 v}{\partial
y_i\partial y_j}+(N-1)r^{\frac{b}{2}-2}\sum_{i=1}^N
x_i\frac{\partial v}{\partial y_i}\nonumber\\
&=&|y|^{-2}\sum^{N}_{i,j=1}y_iy_j\frac{\partial^2 v}{\partial
y_i\partial y_j}+(N-1)|y|^{-2}\sum_{i=1}^N y_i\frac{\partial
v}{\partial y_i}\nonumber\\
&=&\sum_{i,j=1}^N\frac{\partial}{\partial
y_j}\Big(\frac{y_iy_j}{|y|^2}\frac{\partial v}{\partial y_i}\Big).
\end{eqnarray}
Substituting (\ref{kkg88765te}) and $r=|y|^{\frac{2}{2-b}}$ into
(\ref{skkg77tyg6}) results in
\begin{eqnarray}\label{tgkkg7754r55} \triangle_x u(x)
&=&|y|^{-\frac{b(N+2)}{2(2-b)}}\Big(\triangle_y v
+(\frac{b^2}{4}-b)\sum_{i,j=1}^N\frac{\partial}{\partial
y_j}\Big(\frac{y_iy_j}{|y|^2}\frac{\partial v}{\partial y_i}\Big)
-\frac{C_b}{|y|^2}v\Big)\nonumber\\
&=&|y|^{-\frac{b(N+2)}{2(2-b)}}\Big(\sum_{i,j=1}^N\frac{\partial}{\partial
y_j}\Big(A_{ij}(y)\frac{\partial v}{\partial
y_i}\Big)-\frac{C_b}{|y|^2}v\Big).\nonumber
\end{eqnarray}
\hfill$\Box$

Let
\begin{eqnarray}\label{jjgh776ytppl}
H^1_{loc}(\mathbb{R}^N)&=&\{u\ |\ \mbox{for every bounded domain}\
\Omega\subset \mathbb{R}^N,\ \int_{\Omega}|\nabla
u|^2dx+\int_\Omega u^2dx<+\infty\}.
\end{eqnarray}

From the classical Hardy inequality (see, e.g., \cite[Lemma
2.1]{peral}), we deduce that for every bounded $C^1$ domain
$\Omega\subset \mathbb{R}^N$, there exists $C_\Omega>0$ such that,
for every $u\in H^1_{loc}(\mathbb{R}^N)$,
\begin{eqnarray}\label{iig99g8g71q}
\int_\Omega\frac{u^2}{|x|^2}dx\leq C_\Omega(\int_\Omega|\nabla
u|^2dx+\int_\Omega u^2dx)
\end{eqnarray}

\begin{Theorem}\label{kkgn775trff}
If $v\in H^1_{loc}(\mathbb{R}^N)$ is a weak solution of the
equation
\begin{eqnarray}\label{o8867yyttggd}
-\sum_{i,j=1}^N\frac{\partial}{\partial
y_j}\Big(A_{ij}(y)\frac{\partial v}{\partial y_i}\Big)
+V_*v=K_*|v|^{p-2}v\ \mbox{in}\ \mathbb{R}^N,
\end{eqnarray}
i.e., for  every $\psi\in C^\infty_{0}(\mathbb{R}^N),$
\begin{eqnarray}\label{hhfttrgrgd}
\int_{\mathbb{R}^N}\sum^N_{i,j=1}A_{ij}(y)\frac{\partial
v}{\partial y_i}\frac{\partial \psi}{\partial
y_j}dy+\int_{\mathbb{R}^N}V_*(y)v\psi
dy=\int_{\mathbb{R}^N}K_*(y)|v|^{p-2}v\psi dy,
\end{eqnarray}
 and $u$ is
defined by (\ref{lululuggftrr}), then $u\in
H^1_{loc}(\mathbb{R}^N)$ and it is a weak solution of
(\ref{h777arrrrrr}), i.e., for every $\varphi\in
C^\infty_0(\mathbb{R}^N)$,
\begin{eqnarray}\label{oo8rtfgvbbv}
\int_{\mathbb{R}^N}\nabla u\nabla \varphi
dx+\int_{\mathbb{R}^N}V(x)u\varphi
dx=\int_{\mathbb{R}^N}K(x)|u|^{p-2}u\varphi dx.
\end{eqnarray}
\end{Theorem}
\noindent{\bf Proof.} Using the
 spherical coordinates
 \begin{eqnarray}\label{iityh773445}
 &&x_1=r\cos\sigma_1,\nonumber\\
 && x_2=r\sin\sigma_1\cos\sigma_2,\nonumber\\
 &&\ldots\ldots\nonumber\\
&&
x_j=r\sin\sigma_1\sin\sigma_2\cdots\sin\sigma_{j-1}\cos\sigma_j,\
2\leq j\leq N-1,\nonumber\\
 &&\ldots\ldots\nonumber\\
 &&
x_N=r\sin\sigma_1\sin\sigma_2\cdots\sin\sigma_{N-2}\sin\sigma_{N-1},\nonumber
 \end{eqnarray}
where $0\leq\sigma_j<\pi,$ $j=1,2,\ldots, N-2,$
$0\leq\sigma_{N-1}<2\pi,$ we have
\begin{eqnarray}\label{dededeyh}
dx=r^{N-1}f(\sigma)drd\sigma_1\cdots d\sigma_{N-1},\nonumber
\end{eqnarray}
where
$f(\sigma)=\sin^{N-2}\sigma_1\sin^{N-3}\sigma_2\cdots\sin\sigma_{N-2}.$
   Recall that $y=|x|^{-\frac{b}{2}}x.$  Let $R=|y|.$
Then $r=R^{\frac{2}{2-b}}$ and
\begin{eqnarray}\label{mmmmmazx}
dx&=&r^{N-1}f(\sigma)drd\sigma_1\cdots
d\sigma_{N-1}=R^{\frac{2(N-1)}{2-b}}f(\sigma)d(R^{\frac{2}{2-b}})d\sigma_1\cdots
d\sigma_{N-1}\nonumber\\
&=&\frac{2}{2-b}R^{\frac{2N}{2-b}-1}f(\sigma)dRd\sigma_1\cdots
d\sigma_{N-1}=\frac{2}{2-b}|y|^{\frac{bN}{2-b}}dy.
\end{eqnarray}
Here, we used  $dy=R^{N-1}f(\sigma)dRd\sigma_1\cdots
d\sigma_{N-1}$ in the above  last inequality. From
(\ref{kkg77tyg6}), (\ref{mmmmmazx}) and (\ref{iig99g8g71q}), we
deduce that there exists $C>0$ such that for every bounded domain
$\Omega\subset\mathbb{R}^N,$
\begin{eqnarray}\label{9fgbbvn992}
\int_\Omega|\frac{\partial u}{\partial x_i}|^2dx&\leq&
C\int_\Omega r^{-\frac{b(N-2)}{2}-b}\Big(\frac{\partial
v}{\partial
y_i}(|x|^{-b/2}x)\Big)^2dx\nonumber\\
&&+C\int_\Omega r^{-\frac{b(N-2)}{2}-b-4}\Big(x_i\sum_{j=1}^N
x_{j}\frac{\partial
v}{\partial y_j}(|x|^{-b/2}x)\Big)^2dx\nonumber\\
&&+C\int_\Omega
r^{-\frac{b(N-2)}{2}-4}x^2_iv^2(|x|^{-b/2}x)dx\nonumber\\
&=&\frac{2C}{2-b}\int_\Omega \Big(\frac{\partial v(y)}{\partial
y_i}\Big)^2dy+\frac{2C}{2-b}\int_\Omega
\Big(\frac{y_i}{|y|}\sum_{j=1}^N \frac{y_j}{|y|}\frac{\partial
v(y)}{\partial y_j}\Big)^2dy\nonumber\\
&&+\frac{2C}{2-b}\int_\Omega |y|^{-4}y^2_iv^2(y)dy\nonumber\\
&\leq& C''(\int_\Omega |\nabla v|^2dy+\int_\Omega
\frac{v^2}{|y|^2}dy)<+\infty.\nonumber
\end{eqnarray}
Moreover,
\begin{eqnarray}\label{111kkghnv77t6}
\int_\Omega
u^2dx=\int_\Omega|x|^{-\frac{b}{2}(N-2)}v^2(|x|^{-\frac{b}{2}}x)dx=
\frac{2}{2-b}\int_\Omega
|y|^{\frac{2b}{2-b}}v^2(y)dy<+\infty.\nonumber
\end{eqnarray}
Therefore, $u\in H^1_{loc}(\mathbb{R}^N)$. Then, to prove   $u$
satisfies (\ref{oo8rtfgvbbv}) for every $\varphi\in
C^\infty_0(\mathbb{R}^N)$, it suffices to prove that
(\ref{oo8rtfgvbbv}) holds for every $\varphi\in
C^\infty_0(\mathbb{R}^N\setminus\{0\}).$  For $\varphi\in
C^\infty_0(\mathbb{R}^N\setminus\{0\})$, let $\psi\in
C^\infty_0(\mathbb{R}^N\setminus\{0\})$ be such that
\begin{eqnarray}\label{lulbc55ftrr}
\varphi(x)=|x|^{-\frac{b}{4}(N-2)}\psi(|x|^{-\frac{b}{2}}x).\nonumber
\end{eqnarray}
By using the divergence theorem and Lemma \ref{jjfgr665t2w}, we
get that
\begin{eqnarray}\label{00466565y}
&&\int_{\mathbb{R}^N}\nabla u\nabla \varphi
dx\nonumber\\&=&-\int_{\mathbb{R}^N}u\triangle\varphi dx\nonumber\\
&=&-\int_{\mathbb{R}^N}u\cdot|y|^{-\frac{b(N+2)}{2(2-b)}}\Big(\sum_{i,j=1}^N\frac{\partial}{\partial
y_j}\Big(A_{ij}(y)\frac{\partial \psi}{\partial
y_i}\Big)-\frac{C_b}{|y|^2}\psi\Big)dx\nonumber\\
&=&-\int_{\mathbb{R}^N}|x|^{-\frac{b}{4}(N-2)}v(|x|^{-\frac{b}{2}}x)\cdot|y|^{-\frac{b(N+2)}{2(2-b)}}\Big(\sum_{i,j=1}^N\frac{\partial}{\partial
y_j}\Big(A_{ij}(y)\frac{\partial \psi}{\partial
y_i}\Big)-\frac{C_b}{|y|^2}\psi\Big)dx\nonumber\\
&=&-\int_{\mathbb{R}^N}|y|^{-\frac{b(N-2)}{2(2-b)}}v(y)\cdot|y|^{-\frac{b(N+2)}{2(2-b)}}\Big(\sum_{i,j=1}^N\frac{\partial}{\partial
y_j}\Big(A_{ij}(y)\frac{\partial \psi}{\partial
y_i}\Big)-\frac{C_b}{|y|^2}\psi\Big)\frac{2}{2-b}|y|^{\frac{bN}{2-b}}dy\nonumber\\
&=&-\frac{2}{2-b}\int_{\mathbb{R}^N}v\cdot\Big(\sum_{i,j=1}^N\frac{\partial}{\partial
y_j}\Big(A_{ij}(y)\frac{\partial \psi}{\partial
y_i}\Big)-\frac{C_b}{|y|^2}\psi\Big)dy\nonumber\\
&=&\frac{2}{2-b}\int_{\mathbb{R}^N}\sum^N_{i,j=1}A_{ij}(y)\frac{\partial
v}{\partial y_i}\frac{\partial \psi}{\partial
y_j}dy-\frac{2C_b}{2-b}\int_{\mathbb{R}^N}\frac{v\psi}{|y|^2}dy.\nonumber
\end{eqnarray}
Moreover,
\begin{eqnarray}\label{hhfgtbuzhidao}
&&\int_{\mathbb{R}^N}V(x)u\varphi
dx\nonumber\\
&=&\frac{2}{2-b}\int_{\mathbb{R}^N}V(|y|^{\frac{b}{2-b}}y)u(|y|^{\frac{b}{2-b}}y)\varphi(|y|^{\frac{b}{2-b}}y)|y|^{\frac{bN}{2-b}}dy
\nonumber\\
&=&\frac{2}{2-b}\int_{\mathbb{R}^N}|y|^{\frac{2b}{2-b}}V(|y|^{\frac{b}{2-b}}y)\cdot|y|^{\frac{b(N-2)}{2(2-b)}}u(|y|^{\frac{b}{2-b}}y)\cdot
|y|^{\frac{b(N-2)}{2(2-b)}}\varphi(|y|^{\frac{b}{2-b}}y)dy\nonumber\\
&=&\frac{2}{2-b}\int_{\mathbb{R}^N}|y|^{\frac{2b}{2-b}}V(|y|^{\frac{b}{2-b}}y)v(y)\psi(y)dy\nonumber
\end{eqnarray}
and
\begin{eqnarray}\label{jjfh77fythhf1}
&&\int_{\mathbb{R}^N}K(x)|u|^{p-2}u\varphi dx\nonumber\\
&=&\int_{\mathbb{R}^N}K(|y|^{\frac{b}{2-b}}y)|u(|y|^{\frac{b}{2-b}}y)|^{p-2}u(|y|^{\frac{b}{2-b}}y)\varphi(|y|^{\frac{b}{2-b}}y)\frac{2}{2-b}|y|^{\frac{bN}{2-b}}dy\nonumber\\
&=&\frac{2}{2-b}\int_{\mathbb{R}^N}|y|^{\frac{2s}{2-b}}K(|y|^{\frac{b}{2-b}}y)|v(y)|^{p-2}v(y)\psi(y)dy.\nonumber
\end{eqnarray}
Therefore,
\begin{eqnarray}\label{jjfnvnhyg6}
&&\int_{\mathbb{R}^N}\nabla u\nabla \varphi
dx+\int_{\mathbb{R}^N}V(x)u\varphi
dx-\int_{\mathbb{R}^N}K(x)|u|^{p-2}u\varphi dx\nonumber\\
&=&\frac{2}{2-b}\Big(\int_{\mathbb{R}^N}\sum^N_{i,j=1}A_{ij}(y)\frac{\partial
v}{\partial y_i}\frac{\partial \psi}{\partial
y_j}dy-C_b\int_{\mathbb{R}^N}\frac{v\psi}{|y|^2}dy\nonumber\\
&&+\int_{\mathbb{R}^N}|y|^{\frac{2b}{2-b}}V(|y|^{\frac{b}{2-b}}y)v(y)\psi(y)dy\nonumber\\
&&
-\int_{\mathbb{R}^N}|y|^{\frac{2s}{2-b}}K(|y|^{\frac{b}{2-b}}y)|v(y)|^{p-2}v(y)\psi(y)dy\Big)\nonumber\\
&=&\frac{2}{2-b}\Big(\int_{\mathbb{R}^N}\sum^N_{i,j=1}A_{ij}(y)\frac{\partial
v}{\partial y_i}\frac{\partial \psi}{\partial
y_j}dy+\int_{\mathbb{R}^N}V_*(y)v\psi
dy-\int_{\mathbb{R}^N}K_*(y)|v|^{p-2}v\psi
dy\Big)\nonumber\\
&=&0.\nonumber
\end{eqnarray}
This completes the proof. \hfill$\Box$

\medskip

 This theorem implies that the problem of
looking for  solutions of (\ref{h777arrrrrr}) can be reduced  to a
problem of looking for  solutions of (\ref{o8867yyttggd}).

\section{The variational functional for Eq.
(\ref{o8867yyttggd}).}\label{tte54rrer}

 The following inequality is a variant  Hardy inequality.
\begin{lemma}\label{iigy77tyjju} If $ v\in
H^1(\mathbb{R}^N)$, then
\begin{eqnarray}\label{jj88h7tyy}
\int_{\mathbb{R}^N}\frac{|x\cdot\nabla
v|^2}{|x|^2}dx\geq\frac{(N-2)^2}{4}\int_{\mathbb{R}^N}\frac{|
v|^2}{|x|^2}dx.\end{eqnarray}
\end{lemma}
\noindent{\bf Proof.} We only give the proof of (\ref{jj88h7tyy})
for $v\in C^\infty_0(\mathbb{R}^N)$, since $
C^\infty_0(\mathbb{R}^N)$ is dense in $H^1(\mathbb{R}^N).$ For
$v\in C^\infty_0(\mathbb{R}^N)$, we have the following identity
$$|v(x)|^2=-\int^{\infty}_{1}\frac{d}{d\lambda}|v(\lambda x)|^2d\lambda=-2
\int^{\infty}_{1}v(\lambda x)\cdot (x\cdot\nabla v(\lambda
x))d\lambda.$$ By using the H\"older inequality, it follows that
\begin{eqnarray}\label{iit66ghvb}
\int_{\mathbb{R}^N}\frac{|v(x)|^2}{|x|^2}dx&=&-2
\int^{\infty}_{1}\int_{\mathbb{R}^N}\frac{v(\lambda
x)}{|x|^2}\cdot (x\cdot\nabla
v(\lambda x))dxd\lambda \nonumber\\
&=&-2\int^{\infty}_{1}\frac{d\lambda}{\lambda^{N-1}}\int_{\mathbb{R}^N}\frac{v(
x)}{|x|^2}\cdot (x\cdot\nabla v( x))dx\nonumber\\
&=&-\frac{2}{N-2}\int_{\mathbb{R}^N}\frac{v( x)}{|x|^2}\cdot
(x\cdot\nabla v( x))dx\nonumber\\
&\leq&\frac{2}{N-2}(\int_{\mathbb{R}^N}\frac{v^2(
x)}{|x|^2}dx)^{1/2}(\int_{\mathbb{R}^N}\frac{|x\cdot\nabla
v|^2}{|x|^2}dx)^{1/2}.\nonumber
\end{eqnarray}
And then we conclude that
$$\int_{\mathbb{R}^N}\frac{|x\cdot\nabla
v|^2}{|x|^2}dx\geq\frac{(N-2)^2}{4}\int_{\mathbb{R}^N}\frac{|
v|^2}{|x|^2}dx.$$\hfill$\Box$

From the definition of $A_{ij}(x)$ (see (\ref{oo77565ttr})),  it
is easy to verify that, for $u\in H^1(\mathbb{R}^N),$
\begin{eqnarray}\label{mmny66ty}
\int_{\mathbb{R}^N}\sum_{i,j=1}^N A_{ij}(x)\frac{\partial
u}{\partial x_i}\frac{\partial u}{\partial
x_j}dx=\int_{\mathbb{R}^N}|\nabla
u|^2dx+(\frac{b^2}{4}-b)\int_{\mathbb{R}^N}\frac{|x\cdot\nabla
u|^2}{|x|^2}dx.
\end{eqnarray}

\begin{lemma}\label{iidtfhyyrhf}
There exist constants $C_1>0$ and $C_2>0$ such that, for every
$u\in H^1(\mathbb{R}^N)$,
\begin{eqnarray}
C_1||u||^2\leq \int_{\mathbb{R}^N}|\nabla
u|^2dx+(\frac{b^2}{4}-b)\int_{\mathbb{R}^N}\frac{|x\cdot\nabla
u|^2}{|x|^2}dx+\int_{\mathbb{R}^N}V_*(x)|u|^2dx\leq
C_2||u||^2.\nonumber
\end{eqnarray}
\end{lemma}
\noindent{\bf Proof.} From the conditions $\bf (A_1)$ and $\bf
(A_2)$, we  deduce that there exists a constant $C>0$
 such that \begin{eqnarray}\label{uurtf66r5}
 |x|^{\frac{2b}{2-b}}V(|x|^{\frac{b}{2-b}}x)\leq C(1+|x|^{-2}),\ \forall x\in
 \mathbb{R}^N\setminus\{0\}.
 \end{eqnarray}
Since
\begin{eqnarray}\label{iiryyrhfttfg}
\int_{\mathbb{R}^N}V_*(x)|u|^2dx=\int_{\mathbb{R}^N}|x|^{\frac{2b}{2-b}}V(|x|^{\frac{b}{2-b}}x)|u|^2dx+C_b\int_{\mathbb{R}^N}\frac{|
u|^2}{|x|^2}dx,\nonumber
\end{eqnarray}
by (\ref{uurtf66r5}) and the classical Hardy inequality (see,
e.g., \cite{peral})
$$\frac{(N-2)^2}{4}\int_{\mathbb{R}^N}\frac{|
u|^2}{|x|^2}dx\leq\int_{\mathbb{R}^N}|\nabla u|^2dx,\ \forall u\in
H^1(\mathbb{R}^N),$$ we deduce  that there exists a constant $C>0$
 such that
 $$\int_{\mathbb{R}^N}V_*(x)|u|^2dx\leq C||u||^2.$$
This together with the fact that
$\int_{\mathbb{R}^N}\frac{|x\cdot\nabla u|^2}{|x|^2}dx\leq
\int_{\mathbb{R}^N}|\nabla u|^2dx$ yields that there exists a
constant $C_2>0$
 such that
\begin{eqnarray}\label{uufhft6rtf}
&&\int_{\mathbb{R}^N}|\nabla
u|^2dx+(\frac{b^2}{4}-b)\int_{\mathbb{R}^N}\frac{|x\cdot\nabla
u|^2}{|x|^2}dx+\int_{\mathbb{R}^N}V_*(x)|u|^2dx\nonumber\\
&\leq& C_2||u||^2,\ \forall u\in H^1(\mathbb{R}^N).
\end{eqnarray}

If $0<b<2$, then $\frac{b^2}{4}-b<0$ and
\begin{eqnarray}\label{88r6ftfgfg}
\int_{\mathbb{R}^N}|\nabla
u|^2dx+(\frac{b^2}{4}-b)\int_{\mathbb{R}^N}\frac{|x\cdot\nabla
u|^2}{|x|^2}dx&\geq&\int_{\mathbb{R}^N}|\nabla
u|^2dx+(\frac{b^2}{4}-b)\int_{\mathbb{R}^N}|\nabla
u|^2dx\nonumber\\
&=&(1-b/2)^2\int_{\mathbb{R}^N}|\nabla u|^2dx.
\end{eqnarray}
In this case, $C_b=\frac{b}{4}(1-\frac{b}{4})(N-2)^2>0$ and
\begin{eqnarray}\label{iif7f66r992}
\int_{\mathbb{R}^N}V_*(x)|u|^2dx&=&\int_{\mathbb{R}^N}|x|^{\frac{2b}{2-b}}V(|x|^{\frac{b}{2-b}}x)|u|^2dx+C_b\int_{\mathbb{R}^N}\frac{|
u|^2}{|x|^2}dx\nonumber\\
&\geq&\int_{\mathbb{R}^N}|x|^{\frac{2b}{2-b}}V(|x|^{\frac{b}{2-b}}x)|u|^2dx.
\end{eqnarray}
The conditions $\bf (A_1)$ and $\bf (A_2)$ imply that there exists
a constant $C>0$ such that
\begin{eqnarray}\label{gggggggggg}
\int_{\mathbb{R}^N} |\nabla
u|^2dx+\int_{\mathbb{R}^N}|x|^{\frac{2b}{2-b}}V(|x|^{\frac{b}{2-b}}x)u^2dx\geq
C \int_{\mathbb{R}^N} u^2dx .
\end{eqnarray}
Combining  $(\ref{88r6ftfgfg})-(\ref{gggggggggg})$ yields that
there exists a constant $C_1>0$ such that
\begin{eqnarray}\label{k88d6dttfkk9}
&&\int_{\mathbb{R}^N}|\nabla
u|^2dx+(\frac{b^2}{4}-b)\int_{\mathbb{R}^N}\frac{|x\cdot\nabla
u|^2}{|x|^2}dx+\int_{\mathbb{R}^N}V_*(x)|u|^2dx\nonumber\\
&\geq& C_1||u||^2,\ \forall u\in H^1(\mathbb{R}^N).
\end{eqnarray}

If $b<0,$  (\ref{gggggggggg}) still holds. From Lemma
\ref{iigy77tyjju} and (\ref{gggggggggg}), we deduce that there
exists a constant $C_1>0$ such that, for every $u\in
H^1(\mathbb{R}^N)$,
\begin{eqnarray}\label{88r6ftfgfgj}
&&\int_{\mathbb{R}^N}|\nabla
u|^2dx+(\frac{b^2}{4}-b)\int_{\mathbb{R}^N}\frac{|x\cdot\nabla
u|^2}{|x|^2}dx+\int_{\mathbb{R}^N}V_*(x)|u|^2dx\nonumber\\
&=&\int_{\mathbb{R}^N}|\nabla
u|^2dx+(\frac{b^2}{4}-b)\Big(\int_{\mathbb{R}^N}\frac{|x\cdot\nabla
u|^2}{|x|^2}dx-\frac{(N-2)^2}{4}\int_{\mathbb{R}^N}\frac{|
u|^2}{|x|^2}dx\Big)\nonumber\\
&&+\int_{\mathbb{R}^N}|x|^{\frac{2b}{2-b}}V(|x|^{\frac{b}{2-b}}x)|u|^2dx\nonumber\\
&\geq&\int_{\mathbb{R}^N}|\nabla
u|^2dx+\int_{\mathbb{R}^N}|x|^{\frac{2b}{2-b}}V(|x|^{\frac{b}{2-b}}x)|u|^2dx\geq
C_1||u||^2.
\end{eqnarray}
Then the desired result of this lemma follows from
(\ref{uufhft6rtf}), (\ref{k88d6dttfkk9}) and (\ref{88r6ftfgfgj})
immediately.\hfill$\Box$

\medskip

This lemma implies that
\begin{eqnarray}\label{uutyghhgg}
||u||_A=(\int_{\mathbb{R}^N}|\nabla
u|^2dx+(\frac{b^2}{4}-b)\int_{\mathbb{R}^N}\frac{|x\cdot\nabla
u|^2}{|x|^2}dx+\int_{\mathbb{R}^N}V_*(x)|u|^2dx)^{1/2}
\end{eqnarray}
is  equivalent  to the standard  norm  $||\cdot||$ in
$H^1(\mathbb{R}^N)$. We denote the inner product associated with
$||\cdot||_A$ by $(\cdot,\cdot)_A$, i.e.,
\begin{eqnarray}\label{ii9945554}
(u,v)_A
 &=&\int_{\mathbb{R}^N}\nabla u\nabla
vdx+\int_{\mathbb{R}^N}V_*(x)uvdx\nonumber\\
&&+(\frac{b^2}{4}-b)\int_{\mathbb{R}^N}\frac{(x\cdot\nabla
u)(x\cdot\nabla v)}{|x|^2}dx.
\end{eqnarray}By the Sobolev inequality, we have
\begin{eqnarray}\label{uutyrtpi88}
S_A:=\inf_{u\in
H^1(\mathbb{R}^N)\setminus\{0\}}\frac{||u||^2_A}{(\int_{\mathbb{R}^N}|u|^pdx)^{2/p}}>0
\end{eqnarray} and
\begin{eqnarray}\label{iirytfgr1q}
||u||_A\geq S^{\frac{1}{2}}_A(\int_{\mathbb{R}^N}|u|^pdx)^{1/p},\
\forall u\in H^1(\mathbb{R}^N).
\end{eqnarray}
By the condition $\bf(A_1)$ and $\bf(A_2)$, if $0<b<2,$ then $K_*$
is bounded in $\mathbb{R}^N$. Therefore, by (\ref{iirytfgr1q}),
there exists $C>0$ such that
\begin{eqnarray}\label{nnvbf66r5}
(\int_{\mathbb{R}^N}K_*(x)(u^+)^pdx)^{1/p}\leq C||u||_A,\ \forall
u\in H^1(\mathbb{R}^N).
\end{eqnarray}
However, if $b<0$,  $K_*$ has a singularity at $x=0$, i.e.,
\begin{eqnarray}\label{jjvbfgftffr}
K_*(x)\sim |x|^{\frac{2s}{2-b}}K(0),\ \mbox{as}\ |x|\rightarrow 0.
\end{eqnarray}
Recall that $p=2(N-2s/b)/(N-2)$ and $2s/(2-b)>-2s/b$ if $b<0$.
Then by the Hardy-Sobolev inequality (see, for example,
\cite[Lemma 3.2]{ghoussoub}), we deduce that there exists $C>0$
such that (\ref{nnvbf66r5}) still holds. Therefore, the functional
\begin{eqnarray}\label{hhgtrtft33e}
J(u)=\frac{1}{2}||u||^2_A-\frac{1}{p}\int_{\mathbb{R}^N}K_*(x)(u^+)^pdx,\
u\in H^1(\mathbb{R}^N)
\end{eqnarray}
is a $C^2$ functional defined in $H^1(\mathbb{R}^N)$. Moreover, it
easy to check that the Gateaux derivative of $J$ is
$$\langle J'(u), h\rangle=(u,h)_A-\int_{\mathbb{R}^N}K_*(x)(u^+)^{p-1}hdx,\ \forall u,h\in H^1(\mathbb{R}^N)$$
 and  the critical points of
$J$ are nonnegative solutions of (\ref{o8867yyttggd}).

\section{Some minimizing problems}
For $\theta=(\theta_1,\cdots,\theta_N)\in \mathbb{R}^{N}$ with
$|\theta|=1,$  let
\begin{eqnarray}\label{oo775965ttr}
B_{ij}(\theta)=\delta_{ij}+(\frac{b^2}{4}-b)\theta_i\theta_j,\
i,j=1,\cdots,N.
\end{eqnarray}
By this definition, we have, for $u\in H^1(\mathbb{R}^N),$
\begin{eqnarray}\label{kknfgfd6r5r}
\int_{\mathbb{R}^N}\sum_{i,j=1}^NB_{ij}(\theta)\frac{\partial
u}{\partial x_i}\frac{\partial u}{\partial x_j}dx
=\int_{\mathbb{R}^N}|\nabla
u|^2dx+(\frac{b^2}{4}-b)\int_{\mathbb{R}^N}|\theta\cdot\nabla
u|^2dx.
\end{eqnarray}
From
\begin{eqnarray}\label{oofgf66r5tr}
(1+|\frac{b^2}{4}-b|)\int_{\mathbb{R}^N}|\nabla
u|^2dx&\geq&\int_{\mathbb{R}^N}|\nabla
u|^2dx+(\frac{b^2}{4}-b)\int_{\mathbb{R}^N}|\theta\cdot\nabla
u|^2dx\nonumber\\
&\geq&\left\{
\begin{array}{l}
(1-b/2)^2\int_{\mathbb{R}^N}|\nabla u|^2dx, \ 0<b<2\\
\int_{\mathbb{R}^N}|\nabla u|^2dx, \quad b<0,\\
\end{array} \right.\nonumber
\end{eqnarray}we deduce that the norm defined by
\begin{eqnarray}\label{kkfhtyyghg}
||u||_\theta:=(\int_{\mathbb{R}^N}|\nabla
u|^2dx+(\frac{b^2}{4}-b)\int_{\mathbb{R}^N}|\theta\cdot\nabla
u|^2dx +a\int_{\mathbb{R}^N}|u|^2dx)^{1/2}
\end{eqnarray}
is  equivalent to the  standard norm  $||\cdot ||$ in
$H^1(\mathbb{R}^N)$. The inner product corresponding to
$||\cdot||_\theta$ is
$$(u,v)_\theta=\int_{\mathbb{R}^N}\nabla
u\nabla vdx
+a\int_{\mathbb{R}^N}uvdx+(\frac{b^2}{4}-b)\int_{\mathbb{R}^N}(\theta\cdot\nabla
u)(\theta\cdot\nabla v)dx.$$
\begin{lemma}\label{uugfyrttf87}
The infimum
\begin{eqnarray}\label{kkg77t6tt}
\inf_{u\in
H^1(\mathbb{R}^N)\setminus\{0\}}\frac{||u||^2_\theta}{(\int_{\mathbb{R}^N}|u|^pdx)^{2/p}}
\end{eqnarray}
is independent of $\theta\in \mathbb{R}^{N}$ with $|\theta|=1$.
\end{lemma}
\noindent{\bf Proof.} In this proof, we always view a vector in
$\mathbb{R}^N$ as a $1\times N$ matrix. And we use $A^T$ to denote
the conjugate matrix of a matrix $A$.

 For any
$\theta,\theta'\in\mathbb{R}^N$ with $|\theta|=|\theta'|=1$, let
$G$ be an $N\times N$ orthogonal matrix such that $\theta'\cdot
G^T=\theta.$ For any $u\in H^1(\mathbb{R}^N),$ let $v(x)=u(xG),$ $
x\in\mathbb{R}^N.$  The assumption $G$ is an $N\times N$
orthogonal matrix implies that  $GG^T=I,$ where $I$ is the
$N\times N$ identity matrix.  Then it is easy to check that
\begin{eqnarray}\label{kkg77t67yy}
\int_{\mathbb{R}^N}|v|^2dx=\int_{\mathbb{R}^N}|u|^2dx,\
\int_{\mathbb{R}^N}|v|^pdx=\int_{\mathbb{R}^N}|u|^pdx.
\end{eqnarray}
Note that
\begin{eqnarray}\label{1fg66rtglili}
\nabla v(x)=(\nabla u)(xG)\cdot G.
\end{eqnarray}
By  $GG^T=I,$  we have
\begin{eqnarray}\label{iitugyythhg}
|\nabla v(x)|^2&=&\nabla v(x)\cdot (\nabla
v(x))^T\nonumber\\
&=&(\nabla u)(xG)\cdot G\cdot G^T\cdot((\nabla u)(xG))^T=|(\nabla
u)(xG)|^2.\nonumber
\end{eqnarray}
It follows that
\begin{eqnarray}\label{kkghyt776yyy}
\int_{\mathbb{R}^N}|\nabla v(x)|^2dx=\int_{\mathbb{R}^N}|(\nabla
u)(xG)|^2dx=\int_{\mathbb{R}^N}|\nabla u(x)|^2dx.
\end{eqnarray}
By (\ref{1fg66rtglili}) and $\theta'\cdot G^T=\theta,$ we get that
\begin{eqnarray}\label{kkg66gyh}
&&\sum^{N}_{i=1}\theta'_i\frac{\partial v}{\partial x_i}
=\theta'\cdot ((\nabla u)(xG)\cdot G)^T=\theta'\cdot
G^T\cdot((\nabla u)(xG))^T=\theta\cdot((\nabla
u)(xG))^T\nonumber\\
&=&\sum^{N}_{i=1}\theta_i(\frac{\partial u}{\partial
y_i})(xG).\nonumber
\end{eqnarray}
It follows that
\begin{eqnarray}\label{jjghy66ttt}
\int_{\mathbb{R}^N}|\theta'\cdot\nabla
v|^2dx&=&\int_{\mathbb{R}^N}|\sum^{N}_{i=1}\theta'_i\frac{\partial
v}{\partial
x_i}|^2dx\nonumber\\
&=&\int_{\mathbb{R}^N}|\sum^{N}_{i=1}\theta_i(\frac{\partial
u}{\partial y_i})(xG)|^2dx\nonumber\\
&=& \int_{\mathbb{R}^N}|\sum^{N}_{i=1}\theta_i\frac{\partial
u}{\partial x_i}|^2dx= \int_{\mathbb{R}^N}|\theta\cdot\nabla
u|^2dx.
\end{eqnarray}
By (\ref{kkg77t67yy}), (\ref{kkghyt776yyy}) and
(\ref{jjghy66ttt}), we get that $
||v||^2_{\theta'}=||u||^2_\theta. $ This together with
(\ref{kkg77t67yy}) leads to the  result of this lemma.\hfill$\Box$

\medskip

Since  the infimum
 (\ref{kkg77t6tt})
is independent of $\theta\in\mathbb{R}^N$ with $|\theta|=1$, we
 denote it by $S$.
\begin{lemma}\label{fddretrdfdf} Let $S_p$ be the infimum in
 (\ref{olfgcbb5er}). Then $S=(1-b/2)^{\frac{p-2}{p}}S_p$.
\end{lemma}
\noindent{\bf Proof.} Choosing $\theta=(1,0,\cdots,0)$ in
$||\cdot||_\theta$, we have
$$||u||_\theta^2=(1-\frac{b}{2})^2\int_{\mathbb{R}^N}|\frac{\partial u}{\partial x_1}|^2dx+\sum^{N}_{i=2}
\int_{\mathbb{R}^N}|\frac{\partial u}{\partial
x_i}|^2dx+a\int_{\mathbb{R}^N}u^2dx.$$ By Lemma \ref{uugfyrttf87},
we have
\begin{eqnarray}\label{iifdggdv}
S=\inf_{u\in
H^1(\mathbb{R}^N)\setminus\{0\}}\frac{(1-\frac{b}{2})^2\int_{\mathbb{R}^N}|\frac{\partial
u}{\partial x_1}|^2dx+\sum^{N}_{i=2}
\int_{\mathbb{R}^N}|\frac{\partial u}{\partial
x_i}|^2dx+a\int_{\mathbb{R}^N}u^2dx}{(\int_{\mathbb{R}^N}|u|^pdx)^{2/p}}.\nonumber
\end{eqnarray}
Let $$v(x)=u((1-b/2)x_1,x_2,\cdots,x_N),\ x\in \mathbb{R}^N.$$
Then
\begin{eqnarray}\label{uufdterref0}
&&\frac{(1-\frac{b}{2})^2\int_{\mathbb{R}^N}|\frac{\partial
u}{\partial x_1}|^2dx+\sum^{N}_{i=2}
\int_{\mathbb{R}^N}|\frac{\partial u}{\partial
x_i}|^2dx+a\int_{\mathbb{R}^N}u^2dx}{(\int_{\mathbb{R}^N}|u|^pdx)^{2/p}}\nonumber\\
&=&(1-b/2)^{\frac{p-2}{p}}\frac{\int_{\mathbb{R}^N}|\nabla
v|^2dx+a\int_{\mathbb{R}^N}v^2dx}{(\int_{\mathbb{R}^N}|v|^pdx)^{2/p}}.\nonumber
\end{eqnarray}
It follows that
$$S=(1-b/2)^{\frac{p-2}{p}}\inf_{v\in
H^1(\mathbb{R}^N)\setminus\{0\}}\frac{\int_{\mathbb{R}^N}|\nabla
v|^2dx+a\int_{\mathbb{R}^N}v^2dx}{(\int_{\mathbb{R}^N}|v|^pdx)^{2/p}}=(1-b/2)^{\frac{p-2}{p}}S_p.$$
\hfill$\Box$

Since the functionals $||u||^2_\theta$ and
$\int_{\mathbb{R}^N}|u|^pdx$ are invariant by translations, the
same argument as the proof of \cite[Theorem 1.34]{Willem} yields
that there exists a positive minimizer $U_\theta$ for the infimum
$S.$ And from the Lagrange multiplier rule, it is a solution of
\begin{eqnarray}\label{hhftr6yyf} -\sum_{i,j=1}^N\frac{\partial}{\partial
y_j}\Big(B_{ij}(\theta)\frac{\partial u}{\partial
y_i}\Big)+au=S(u^+)^{p-1}\ \mbox{in}\ \mathbb{R}^N\nonumber
\end{eqnarray}
and $(\mu/S)^{-1/(p-2)}U_\theta$ is a solution of
 \begin{eqnarray}\label{hhftr6yyf}
-\sum_{i,j=1}^N\frac{\partial}{\partial
y_j}\Big(B_{ij}(\theta)\frac{\partial u}{\partial
y_i}\Big)+au=\mu(u^+)^{p-1}\ \mbox{in}\ \mathbb{R}^N.
\end{eqnarray}
In the next section, we shall show that  Eq.(\ref{hhftr6yyf}) is
the "limit" equation of
\begin{eqnarray}\label{iiijkmmn}
-\sum_{i,j=1}^N\frac{\partial}{\partial
y_j}\Big(A_{ij}(x)\frac{\partial u}{\partial
y_i}\Big)+V_*(x)u=K_*(x)(u^+)^{p-1}\ \mbox{in}\ \mathbb{R}^N.
\end{eqnarray}
 It is easy to verify that
\begin{eqnarray}\label{uugf6tyyg11}
J_{\theta}(u)
=\frac{1}{2}||u||^2_\theta-\frac{\mu}{p}\int_{\mathbb{R}^N}(u^+)^pdx,\
u\in H^1(\mathbb{R}^N),
\end{eqnarray}
 is a $C^2$ functional defined in $H^1(\mathbb{R}^N)$, the Gateaux derivative of $J_\theta$ is
$$\langle J'_\theta(u), h\rangle=(u,h)_\theta-\mu\int_{\mathbb{R}^N}(u^+)^{p-1}hdx,\ \forall u,h\in H^1(\mathbb{R}^N).$$
 and the
 critical points of this functional
are  solutions of (\ref{hhftr6yyf}).

\begin{lemma}\label{iighjgytyf}
Let $\theta\in\mathbb{R}^N$  satisfy $|\theta|=1.$ If $u\neq 0$ is
a critical point of $J_\theta$, then
\begin{eqnarray}\label{kkgh77tytyg}
J_\theta(u)\geq
(\frac{1}{2}-\frac{1}{p})\mu^{-\frac{2}{p-2}}S^{\frac{p}{p-2}}.
\end{eqnarray}
\end{lemma}
\noindent{\bf Proof.} Since $u $ is a critical point of
$J_\theta$, we have
\begin{eqnarray}\label{llmaxxtggf}
&&0=\langle J'_\theta(u),u\rangle
=||u||^2_\theta-\mu\int_{\mathbb{R}^N}(u^+)^pdx.
\end{eqnarray}
It follows that \begin{eqnarray}\label{km77ytg00}
J_\theta(u)=(\frac{1}{2}-\frac{1}{p})\mu\int_{\mathbb{R}^N}(u^+)^pdx.
\end{eqnarray}
Since $u\neq 0,$ by
$||u||^2_\theta=\mu\int_{\mathbb{R}^N}(u^+)^pdx$ and
$||u||^2_\theta\geq S(\int_{\mathbb{R}^N}(u^+)^pdx)^{2/p}$, we get
that $$\int_{\mathbb{R}^N}(u^+)^pdx\geq(S/\mu)^{p/(p-2)}.$$ This
together with (\ref{km77ytg00}) yields the  result of this
lemma.\hfill$\Box$

\section{ Palais-Smale conditions for the functional $J$.}

Recall that $J$ is the functional defined by (\ref{hhgtrtft33e}).
By a $(PS)_c$ sequence of $J,$  we mean a sequence $\{u_n\}\subset
H^1(\mathbb{R}^N)$ such that $J(u_n)\rightarrow c$ and
$J'(u_n)\rightarrow 0$ in $H^{-1}(\mathbb{R}^N)$ as
$n\rightarrow\infty,$  where $H^{-1}(\mathbb{R}^{N})$ denotes the
dual space of $H^{1}(\mathbb{R}^{N})$. $J$ is called satisfying
$(PS)_c$ condition if every  $(PS)_c$ sequence of $J$ contains a
convergent subsequence in $H^1(\mathbb{R}^N)$.

 Our main result in this section  reads as follows:
\begin{Theorem}\label{uut6ghvbbo}
Under the assumptions  $\bf (A_1)$ and $\bf (A_2)$, let
$\{u_n\}\subset H^1(\mathbb{R}^N)$ be a $(PS)_c$ sequence of $J$.
Then replacing $\{u_n\}$ if necessary by a subsequence, there
exist a solution $u_0\in H^1(\mathbb{R}^N)$ of
Eq.(\ref{iiijkmmn}), a finite sequence $\{\theta_l\in\mathbb{R}^N\
|\ |\theta_l|=1,\ 1\leq l\leq k\}$, $k$ functions $\{u_l\ |\ 1\leq
i\leq k\}\subset H^1(\mathbb{R}^N)$
 and $k$ sequences $\{y^l_n\}\subset\mathbb{R}^N$ satisfying:
 \begin{itemize}
\item[{(i).}] $-\sum_{i,j=1}^N\frac{\partial}{\partial
y_j}\Big(B_{ij}(\theta_l)\frac{\partial u_l}{\partial
y_i}\Big)+au_l=\mu(u^+_l)^{p-1}\ \mbox{in}\ \mathbb{R}^N,$
\item[{(ii).}] $|y^l_n|\rightarrow\infty,\
|y^l_{n}-y^{l'}_n|\rightarrow\infty,\ l\neq l',\
n\rightarrow\infty,$\item[{(iii).}]$
||u_n-u_0-\sum^{k}_{l=1}u_l(\cdot-y^l_n)||\rightarrow 0,$
\item[{(iv).}] $J(u_0)+\sum^{l}_{i=1}J_{\theta_l}(u_l)=c.$
\end{itemize}
\end{Theorem}

This theorem gives a precise representation of $(PS)_c$ sequence
for the functional $J$. Through it, partial compactness for $J$
can be regained (see Corollary \ref{kkg77t6yy}).

To prove this theorem, we need some lemmas.  Our proof of this
theorem is inspired by the proof of \cite[Theorem 8.4]{Willem}.

\begin{lemma}\label{nnvbf5rtf}
Let $u\in H^1(\mathbb{R}^N)$. Then for any sequence
$\{y_n\}\subset\mathbb{R}^N,$
$$
\lim_{R\rightarrow\infty}\sup_{n}\int_{|x|>R}K_*(x+y_n)|u|^pdx=0.$$
If $|y_n|\rightarrow\infty,$ $n\rightarrow\infty$, then $$
\lim_{n\rightarrow\infty}\int_{\mathbb{R}^N}|K_*(x+y_n)-\mu|\cdot|u|^pdx=0.$$
\end{lemma}
\noindent{\bf Proof.} If $2>b>0,$ then $K_*$ is bounded in
$\mathbb{R}^N$. In this case, the result of this lemma is obvious.
If $b<0$, then $K_*(x)\sim |x|^{\frac{2s}{2-b}}K(0)$ as
$|x|\rightarrow 0.$ Since $2s/(2-b)>-2s/b$,  by Lemma 3.2 of
\cite{ghoussoub}, the map $v\mapsto K^{1/p}_*v$ from
$H^1(\mathbb{R}^N)\rightarrow L^p_{loc}(\mathbb{R}^N)$ is compact.
Therefore, for any $\epsilon>0$, there exists $\delta_\epsilon>0$
such that
\begin{eqnarray}\label{nvnvbtfll}
\sup_{n}\int_{|x|\leq\delta_\epsilon}K_*(x)|u(x-y_n)|^pdx\leq\epsilon.\nonumber
\end{eqnarray}
And there exists $D(\epsilon)>0$ depending only on $\epsilon$ such
that $K_*(x)\leq D(\epsilon),\ |x|\geq\delta_\epsilon.$ Then for
every $n,$
\begin{eqnarray}
&&\int_{|x|>R}K_*(x+y_n)|u|^pdx\nonumber\\
&&\leq\int_{\{x\ |\ |x+y_n|\leq\delta_\epsilon,\
|x|>R\}}K_*(x+y_n)|u|^pdx+\int_{\{x\ |\ |x+y_n|>\delta_\epsilon,\
|x|>R\}}K_*(x+y_n)|u|^pdx\nonumber\\
&\leq&\epsilon+C(\epsilon)\int_{ |x|>R}|u|^pdx.\nonumber
\end{eqnarray}
It follows that $
\limsup_{R\rightarrow\infty}\sup_{n}\int_{|x|>R}K_*(x+y_n)|u|^pdx\leq\epsilon.$
Now let $\epsilon\rightarrow 0.$

Using the same argument in the above, for any $\epsilon>0,$ there
exist $\delta_\epsilon$ and $D(\epsilon)$ such that
\begin{eqnarray}\label{nvnvbtfll}
\sup_{n}\int_{|x+y_n|\leq\delta_\epsilon}|K_*(x+y_n)-\mu|\cdot|u|^pdx\leq\epsilon\nonumber
\end{eqnarray}
and
$$|K_*(x+y_n)-\mu|\cdot|u|^pdx\leq (D(\epsilon)+\mu)|u|^p, \ |x+y_n|\geq\delta_\epsilon.$$
Since $y_n\rightarrow\infty$, we have  $\lim K_*(x+y_n)=\mu$. Then
using the Lebesgue theorem and the above two inequalities, we get
that
$$
\limsup_{n\rightarrow\infty}\int_{\mathbb{R}^N}|K_*(x+y_n)-\mu|\cdot|u|^pdx\leq\epsilon.$$
Let $\epsilon\rightarrow 0$. Then we get the desired result of
this lemma.
 \hfill$\Box$

  \begin{lemma}\label{llg66ghyhhv}
 Let $\rho>0$. If $\{u_{n}\}$ is bounded in
$H^{1}(\mathbb{R}^{N})$ and if
\begin{equation}\label{nncb5f5f11}
\sup_{y\in \mathbb{R}^{N}}\int_{B(y,
\rho)}|u_{n}|^{2}dx\rightarrow 0, \ n\rightarrow\infty,
\end{equation}
then $k^{1/p}_*u_{n}\rightarrow 0$ in $L^{p}(\mathbb{R}^{N})$.
\end{lemma}

\noindent{\bf Proof.} Since $2s/(2-b)>-2s/b$,  by Lemma 3.2 of
\cite{ghoussoub}, the map $v\mapsto K^{1/p}_*v$ from
$H^1(\mathbb{R}^N)\rightarrow L^p_{loc}(\mathbb{R}^N)$ is compact.
Therefore, for any $\epsilon>0$, there exists $\delta_\epsilon>0$
such that
\begin{eqnarray}\label{nvnvbtfll}
\sup_{n}\int_{|x|\leq\delta_\epsilon}K_*(x)|u_n|^pdx\leq\epsilon.\nonumber
\end{eqnarray}
And  there exists $D(\epsilon)>0$ depending only on $\epsilon$
such that $K_*(x)\leq D(\epsilon),\ |x|\geq\delta_\epsilon.$ By
(\ref{nncb5f5f11}) and the Lions Lemma (see, for example,
\cite[Lemma 1.21]{Willem}), we get that
\begin{eqnarray}\label{hhcgfffdf}
\int_{|x|\geq\delta_\epsilon}K_*(x)|u_n|^pdx \leq
D(\epsilon)\int_{\mathbb{R}^N}|u_n|^pdx\rightarrow 0,\
n\rightarrow\infty.\nonumber
\end{eqnarray}
Therefore,$\limsup_{n\rightarrow\infty}\int_{\mathbb{R}^N}K_*(x)|u_n|^pdx\leq\epsilon.$
Now let $\epsilon\rightarrow 0.$ \hfill$\Box$

\begin{lemma}\label{ooh8y77yuyy} Let
$\{y_n\}\subset\mathbb{R}^N.$
 If $u_{n}\rightharpoonup u$ in
$H^{1}(\mathbb{R}^{N})$, then
$$K_*(x+y_n)(u_{n}^+)^{p-1}-K_*(x+y_n)((u_{n}-u)^+)^{p-1}
-K_*(x+y_n)(u^+)^{p-1}\rightarrow0 \ \mbox{in} \
H^{-1}(\mathbb{R}^{N}).$$
\end{lemma}
One can follow the proof of  \cite[Lemma 8.1]{Willem} step by step
and use Lemma \ref{nnvbf5rtf} to give the  proof of this lemma.

 \medskip

The following Lemma is a variant Br\'ezis-Lieb Lemma (see
\cite{BLI}) and its proof is similar to that of \cite[Lemma 1.32
]{Willem}.
\begin{lemma}\label{honhtton} Let
$\{u_{n}\}\subset H^{1}(\mathbb{R}^{N})$ and
$\{y_n\}\subset\mathbb{R}^N.$
 If \begin{description}
 \item{a)} $\{u_{n}\}$ is bounded in $ H^{1}(\mathbb{R}^{N}),$
 \item{b)}  $u_{n}\rightarrow u$ a.e. on $\mathbb{R}^{N},$
  then
\end{description}
$$\displaystyle\lim_{n\rightarrow \infty}
\int_{\mathbb{R}^{N}}K_*(x+y_n)\cdot|(u^+_n)^{p}-((u_n-u)^+)^{p}-(u^+)^{p}|dx=0.$$
\end{lemma}
\noindent{\bf Proof.} Let $j(t)=\left\{\begin{array}{l}t^p, \ t\geq 0\\
0,\ \ t<0 \\ \end{array} \right..$ Then $j$ is  a convex function.
From \cite[Lemma 3 ]{BLI}, we have for any $\epsilon>0,$ there
exists $C(\epsilon)>0$ such that for all $a,b\in\mathbb{R},$
\begin{eqnarray}\label{nnvb6ftrg}
|j(a+b)-j(b)|\leq\epsilon j(a)+C(\epsilon)j(b).
\end{eqnarray}
Hence
\begin{eqnarray}
f_n^\epsilon&:=&\Big(K_*(x+y_n)\cdot|(u^+_n)^{p}-((u_n-u)^+)^{p}-(u^+)^{p}|-\epsilon
K_*(x+y_n)\cdot((u_n-u)^{+})^p\Big)^+\nonumber\\
&\leq& (1+C(\epsilon))K_*(x+y_n)\cdot(u^+)^p.\nonumber
\end{eqnarray}
By   Lemma 3.2 of \cite{ghoussoub}, the map $v\mapsto K^{1/p}_*v$
from $H^1(\mathbb{R}^N)\rightarrow L^p_{loc}(\mathbb{R}^N)$ is
compact. We get that there exists $\delta_\epsilon>0$ such that
for any $n,$
\begin{eqnarray}\label{hhfgvcrf}
\int_{|x+y_n|<\delta_\epsilon}f^\epsilon_ndx<\epsilon.
\end{eqnarray}
And there exists $D(\epsilon)>0$ depending only on $\epsilon$ such
that $K_*(x)\leq D(\epsilon),\ |x|\geq\delta_\epsilon.$ Then
\begin{eqnarray}\label{hhfgvb66frtfg}
f_n^\epsilon\leq (1+C(\epsilon))D(\epsilon)\cdot(u^+)^p,\
|x+y_n|\geq\delta_\epsilon.\nonumber
\end{eqnarray}
By the Lebesgue theorem,
$\int_{|x+y_n|\geq\delta_\epsilon}f_n^\epsilon dx\rightarrow 0,\
n\rightarrow\infty.$ This together with (\ref{hhfgvcrf}) yields $$
\limsup_{n\rightarrow\infty}\int_{\mathbb{R}^N}f_n^\epsilon
dx\leq\epsilon.$$ The left proof is the same as the proof of
\cite[Lemma 1.32]{Willem}.\hfill$\Box$

\begin{lemma}\label{iig8uuhn} If \begin{eqnarray}
& &u_{n}\rightharpoonup u  \ \mbox{in} \  H^{1}(\mathbb{R}^N),\nonumber \\
& &u_{n}\rightarrow u  \  \mbox{a.e. on}  \ \mathbb{R}^N, \nonumber \\
& &J(u_{n})\rightarrow c, \nonumber \\
& &J'(u_{n})\rightarrow 0  \ \mbox{in} \  H^{-1}(\mathbb{R}^N),
\nonumber
\end{eqnarray}
then $J'(u)=0$ in $H^{-1}(\mathbb{R}^N)$ and $v_{n}:= u_{n}-u$ is
such that \begin{eqnarray}
& &||v_{n}||^{2}_A=||u_{n}||^{2}_A-||u||^{2}_A+o(1),\nonumber \\
& &J(v_{n})\rightarrow c-J(u),\nonumber \\
& &J'(v_{n})\rightarrow 0  \ \mbox{in} \ H^{-1}(\mathbb{R}^N).
\nonumber
\end{eqnarray}
\end{lemma}

\noindent{\bf Proof.} 1). Since $u_n\rightharpoonup u$ in
$H^1(\mathbb{R}^N)$, we get that, as $n\rightarrow\infty,$
\begin{eqnarray}\label{jguyyt0o}
||v_n||^2_A-||u_n||^2_A=(u_n-u,u_n-u)_A-||u_n||^2_A=||u||^2_A-2(u_n,u)_A\rightarrow-||u||^2_A.\nonumber
\end{eqnarray}
Therefore,
\begin{eqnarray}\label{0oortbv}
||v_{n}||^{2}_A=||u_{n}||^{2}_A-||u||^{2}_A+o(1).
\end{eqnarray}

2).  Lemma \ref{honhtton} implies
\begin{eqnarray}\label{og77tygh}
\int_{\mathbb{R}^N}K_*(x)(v^+_n)^pdx=\int_{\mathbb{R}^N}K_*(x)(u^+_n)^pdx-\int_{\mathbb{R}^N}K_*(x)(u^+)^pdx+o(1).
\end{eqnarray}
By (\ref{0oortbv}), (\ref{og77tygh}) and the assumption
$J(u_n)\rightarrow c$, we get that
$$J(v_{n})\rightarrow c-J(u),\ n\rightarrow\infty.$$

3). Since $J'(u_n)\rightarrow 0$ in $H^{-1}(\mathbb{R}^N)$ and
$u_n\rightharpoonup u,$ it is easy to verify that $J'(u)=0.$ For
$h\in H^1(\mathbb{R}^N)$,
\begin{eqnarray}\label{nbuuhjhj}
\langle J'(v_n),
h\rangle&=&(v_n,h)_A-\int_{\mathbb{R}^N}K_*(x)(v^+_n)^{p-1}hdx\nonumber\\
&=&(u_n, h)_A-(u,h)_A-\int_{\mathbb{R}^N}K_*(x)(v^+_n)^{p-1}hdx.
\end{eqnarray}
By Lemma \ref{ooh8y77yuyy}, we have
\begin{eqnarray}\label{kkhnb88t6ty}
&&\sup_{||h||\leq
1}|\int_{\mathbb{R}^N}K_*(x)(v^+_n)^{p-1}hdx-\int_{\mathbb{R}^N}K_*(x)(u^+_n)^{p-1}hdx+\int_{\mathbb{R}^N}K_*(x)(u^+)^{p-1}hdx|\nonumber\\
&&\rightarrow 0,\ n\rightarrow\infty.
\end{eqnarray}
Combining (\ref{nbuuhjhj}) and  (\ref{kkhnb88t6ty}) leads to
$J'(v_n)=J'(u_n)-J'(u)+o(1).$ Then by  $J'(u_n)\rightarrow 0$ in
$H^{-1}(\mathbb{R}^N)$ and  $J'(u)=0$, we obtain that
$J'(v_n)\rightarrow 0$ in $H^{-1}(\mathbb{R}^N)$.
 \hfill$\Box$

\begin{lemma}\label{dondgdonguu}
If $|y_n|\rightarrow\infty$ and as $n\rightarrow\infty,$
\begin{eqnarray}&&u_n(\cdot+y_n)\rightharpoonup u\ \mbox{in}\
H^1(\mathbb{R}^N),\nonumber\\
&&u_n(\cdot+y_n)\rightarrow u\ a.e. \ \mbox{on}\
\mathbb{R}^N,\nonumber\\
&&J(u_n)\rightarrow c,\nonumber\\
&&J'(u_n)\rightarrow 0\ \mbox{in}\ H^{-1}(\mathbb{R}^N),\nonumber
\end{eqnarray}
then there exists $\theta\in\mathbb{R}^N$ with $|\theta|=1$ such
that $J'_\theta(u)=0$ and $v_n=u_n-u(\cdot-y_n)$ is such that
 \begin{eqnarray}
& &||v_{n}||^{2}=||u_{n}||^{2}-||u||^{2}+o(1),\nonumber \\
& &J(v_{n})\rightarrow c-J_\theta(u),\nonumber \\
& &J'(v_{n})\rightarrow 0  \ \mbox{in} \ H^{-1}(\mathbb{R}^N).
\nonumber
\end{eqnarray}
\end{lemma}
\noindent{\bf Proof.} We divide the proof into several steps.

 1). Since $u_n(\cdot+y_n)\rightharpoonup u$
in $ H^1(\mathbb{R}^N)$, it is clear that
$$||v_n||^2=||v_n(\cdot+y_n)||^2=||u_n(\cdot+y_n)||^2+||u||^2-2(u_n(\cdot+y_n),u)=||u_n||^2-||u||^2+o(1).$$

2). For any $h\in H^1(\mathbb{R}^N)$,
\begin{eqnarray}\label{hftryfgggf}
\langle J'(u_n),
h(\cdot-y_n)\rangle=(u_n,h(\cdot-y_n))_A-\int_{\mathbb{R}^N}K_*(x)(u_n^+)^{p-1}h(\cdot-y_n)dx.
\end{eqnarray}
By the definition of the inner product $(\cdot,\cdot)_A$ (see
(\ref{ii9945554})), we have
\begin{eqnarray}\label{ii9945554es}
&&(u_n,h(\cdot-y_n))_A\nonumber\\
&=&\int_{\mathbb{R}^N}\nabla u_n\nabla
h(\cdot-y_n)dx+(\frac{b^2}{4}-b)\int_{\mathbb{R}^N}\frac{(x\cdot\nabla
u_n)(x\cdot\nabla
h(\cdot-y_n))}{|x|^2}dx\nonumber\\
&&+\int_{\mathbb{R}^N}V_*(x)u_nh(\cdot-y_n)dx\nonumber\\
&=&\int_{\mathbb{R}^N}\nabla u_n(\cdot+y_n)\nabla
hdx+a\int_{\mathbb{R}^N}
u_n(\cdot+y_n)\cdot hdx\nonumber\\
&&+\int_{\mathbb{R}^N}(V_*(x+y_n)-a)u_n(\cdot+y_n)\cdot
hdx\nonumber\\
&&+(\frac{b^2}{4}-b)\int_{\mathbb{R}^N}\Big(\frac{\frac{x}{|y_n|}+\frac{y_n}{|y_n|}}{|\frac{x}{|y_n|}+\frac{y_n}{|y_n|}|}\cdot\nabla
u_n(\cdot+y_n)\Big)\Big(\frac{\frac{x}{|y_n|}+\frac{y_n}{|y_n|}}{|\frac{x}{|y_n|}+\frac{y_n}{|y_n|}|}\cdot\nabla
h\Big)dx\nonumber\\
&:=&I+II+III.
\end{eqnarray}

Since $u_n(\cdot+y_n)\rightharpoonup u$ in $H^1(\mathbb{R}^N)$, we
have
\begin{eqnarray}\label{kkfttrg6645}
&&I=\int_{\mathbb{R}^N}\nabla u_n(\cdot+y_n)\nabla
hdx+a\int_{\mathbb{R}^N} u_n(\cdot+y_n)\cdot
hdx=(u_n(\cdot+y_n),h)\nonumber\\
&&\rightarrow  \int_{\mathbb{R}^N}\nabla u\nabla
hdx+a\int_{\mathbb{R}^N} uhdx,\ n\rightarrow\infty.
\end{eqnarray}

By the assumption $\bf (A_2)$ and the definition of $V_*$, we have
$\lim_{|x|\rightarrow\infty}V_*(x)=a$. This yields
\begin{eqnarray}\label{mmbnvyyft}
\sup_{n}\int_{|x|\geq R}|V_*(x)-a|\cdot|h(x-y_n)|^2dx\rightarrow
0,\ R\rightarrow\infty.\nonumber
\end{eqnarray}
Moreover, together with (\ref{iig99g8g71q}) and the fact that
$|y_n|\rightarrow\infty$ yields that for any fixed $R>0$
\begin{eqnarray}\label{uuthrffegbbvb}
&&\int_{|x|<R}|V_*(x)-a|\cdot|h(\cdot-y_n)|^2dx\nonumber\\
&\leq& C(\int_{|x|<R}|\nabla
h(\cdot-y_n)|^2dx+\int_{|x|<R}|h(\cdot-y_n)|^2dx)\rightarrow0,\
n\rightarrow\infty.\nonumber
\end{eqnarray}
Combining the above two limits leads to
\begin{eqnarray}\label{uuthgbbvb}
&&\int_{\mathbb{R}^N}|V_*(x+y_n)-a|\cdot|h|^2dx\rightarrow0,\
n\rightarrow\infty.
\end{eqnarray}

 By (\ref{uuthgbbvb}) and
the H\"older inequality, we have
\begin{eqnarray}\label{kkkmnhyyt}
&&|II|=|\int_{\mathbb{R}^N}(V_*(x+y_n)-a)u_n(\cdot+y_n)\cdot
hdx|\nonumber\\
&\leq&(\int_{\mathbb{R}^N}|V_*(x+y_n)-a|u^2_n(\cdot+y_n)dx)^{\frac{1}{2}}(\int_{\mathbb{R}^N}|V_*(x+y_n)-a|h^2dx)^{\frac{1}{2}}\nonumber\\
&\leq&
C(\int_{\mathbb{R}^N}|V_*(x+y_n)-a|h^2dx)^{\frac{1}{2}}\rightarrow
0,\ n\rightarrow\infty.
\end{eqnarray}

Since $\nabla h\in L^2(\mathbb{R}^N)$, for any $\epsilon>0,$ there
exists $R_\epsilon>0$ such that
\begin{eqnarray}\label{uurtrggfgf}
\int_{\mathbb{R}^N\setminus \{|x|<R_\epsilon\}}|\nabla
h|^2dx<\epsilon.\nonumber
\end{eqnarray}
It follows that
\begin{eqnarray}\label{jjjjjhuyuyu}
\int_{\mathbb{R}^N\setminus
\{|x|<R_\epsilon\}}\frac{|(\frac{x}{|y_n|}+\frac{y_n}{|y_n|})\cdot\nabla
h|^2}{|\frac{x}{|y_n|}+\frac{y_n}{|y_n|}|^2}dx\leq\int_{\mathbb{R}^N\setminus
\{|x|<R_\epsilon\}}|\nabla h|^2dx<\epsilon.
\end{eqnarray}
Then
\begin{eqnarray}\label{iithhgnvb766}
&&|\int_{\mathbb{R}^N\setminus\{|x|<
R_\epsilon\}}\Big(\frac{\frac{x}{|y_n|}+\frac{y_n}{|y_n|}}{|\frac{x}{|y_n|}+\frac{y_n}{|y_n|}|}\cdot\nabla
u_n(\cdot+y_n)\Big)\Big(\frac{\frac{x}{|y_n|}+\frac{y_n}{|y_n|}}{|\frac{x}{|y_n|}+\frac{y_n}{|y_n|}|}\cdot\nabla
h\Big)dx|\nonumber\\
&\leq&(\int_{\mathbb{R}^N\setminus
\{|x|<R_\epsilon\}}\frac{|(\frac{x}{|y_n|}+\frac{y_n}{|y_n|})\cdot\nabla
 u_n(\cdot+y_n)|^2}{|\frac{x}{|y_n|}+\frac{y_n}{|y_n|}|^2}dx)^{\frac{1}{2}}\nonumber\\
&&\times(\int_{\mathbb{R}^N\setminus
\{|x|<R_\epsilon\}}\frac{|(\frac{x}{|y_n|}+\frac{y_n}{|y_n|})\cdot\nabla
h|^2}{|\frac{x}{|y_n|}+\frac{y_n}{|y_n|}|^2}dx)^{\frac{1}{2}}\nonumber\\
&\leq& (\int_{\mathbb{R}^N}|\nabla
u_n|^2dx)^{1/2}(\int_{\mathbb{R}^N\setminus\{|x|<R_\epsilon\}}|\nabla
h|^2dx)^{1/2}\leq C\epsilon
\end{eqnarray}
 where the constant $C$ is independent of $\epsilon$ and $n.$
There exist a subsequence of $y_n/|y_n|$, denoted by itself for
convenience, and $\theta\in \mathbb{R}^N$ with $|\theta|=1$ such
that $y_n/|y_n|\rightarrow \theta$ as $n\rightarrow\infty.$ Then
by $|y_n|\rightarrow\infty$, we get that, as $n\rightarrow\infty$,
$$\frac{x}{|y_n|}+\frac{y_n}{|y_n|}\rightarrow\theta,\ a.e.\ \mbox{on} \ \mathbb{R}^N$$
and $\frac{x}{|y_n|}+\frac{y_n}{|y_n|}$ converges to $\theta$
uniformly for $|x|< R_\epsilon.$ Therefore, there exists
$N_\epsilon$ such that, when $n>N_\epsilon$,
\begin{eqnarray}\label{tgyyhhyy}
&&|\int_{\{|x|<
R_\epsilon\}}\Big(\frac{\frac{x}{|y_n|}+\frac{y_n}{|y_n|}}{|\frac{x}{|y_n|}+\frac{y_n}{|y_n|}|}\cdot\nabla
u_n(\cdot+y_n)\Big)\Big(\frac{\frac{x}{|y_n|}+\frac{y_n}{|y_n|}}{|\frac{x}{|y_n|}+\frac{y_n}{|y_n|}|}\cdot\nabla
h\Big)dx\nonumber\\
&&\quad-\int_{ \{|x|<R_\epsilon\}}(\theta\cdot\nabla
u_n(\cdot+y_n))(\theta\cdot\nabla h)dx|<\epsilon.
\end{eqnarray}
Since $u_n(\cdot+y_n)\rightharpoonup u$ in $H^1(\mathbb{R}^N)$, we
have $\nabla u_n(\cdot+y_n)\rightharpoonup \nabla u$ in
$L^2(\mathbb{R}^N)$. It implies that
$$\int_{ \{|x|<R_\epsilon\}}(\theta\cdot\nabla
u_n(\cdot+y_n))(\theta\cdot\nabla h)dx\rightarrow\int_{
\{|x|<R_\epsilon\}}(\theta\cdot\nabla u)(\theta\cdot\nabla h)dx,\
n\rightarrow\infty.$$ This together with (\ref{iithhgnvb766}),
$(\ref{tgyyhhyy})$ and
\begin{eqnarray}\label{iiiierdf} \int_{
\mathbb{R}^N\setminus\{|x|<R_\epsilon\}}|\theta\cdot\nabla
h|^2dx\leq\int_{\mathbb{R}^N\setminus \{|x|<R_\epsilon\}}|\nabla
h|^2dx<\epsilon,\nonumber
\end{eqnarray}
 yields that there
exists $N'_\epsilon>0$ such that when $n>N'_\epsilon$,
\begin{eqnarray}\label{shandongs}
|\int_{\mathbb{R}^N}\Big(\frac{\frac{x}{|y_n|}+\frac{y_n}{|y_n|}}{|\frac{x}{|y_n|}+\frac{y_n}{|y_n|}|}\cdot\nabla
u_n(\cdot+y_n)\Big)\Big(\frac{\frac{x}{|y_n|}+\frac{y_n}{|y_n|}}{|\frac{x}{|y_n|}+\frac{y_n}{|y_n|}|}\cdot\nabla
h\Big)dx -\int_{ \mathbb{R}^N}(\theta\cdot\nabla
u)(\theta\cdot\nabla h)dx|<(4+C)\epsilon.\nonumber
\end{eqnarray}
Thus
\begin{eqnarray}\label{mmgyrttr6}
III\rightarrow (\frac{b^2}{4}-b)\int_{
\mathbb{R}^N}(\theta\cdot\nabla u)(\theta\cdot\nabla h)dx,\
n\rightarrow\infty.
\end{eqnarray}
Combining  (\ref{kkfttrg6645}), (\ref{kkkmnhyyt}) and
(\ref{mmgyrttr6}) leads to
\begin{eqnarray}\label{mxxxcxcxx}
&&(u_n,h(\cdot-y_n))_A\nonumber\\
&=&\int_{\mathbb{R}^N}\nabla u\nabla hdx+a\int_{\mathbb{R}^N}
uhdx+(\frac{b^2}{4}-b)\int_{ \mathbb{R}^N}(\theta\cdot\nabla
u)(\theta\cdot\nabla h)dx+o(1)\nonumber\\
&=&(u,h)_\theta+o(1).
\end{eqnarray}

We obtain by the H\"older inequality and Lemma \ref{nnvbf5rtf}
that, as $n\rightarrow\infty,$
\begin{eqnarray}\label{mmbnbvx}
&&|\int_{\mathbb{R}^N}K_*(x+y_n)(u^+_n(\cdot+y_n))^{p-1}hdx-\mu\int_{\mathbb{R}^N}(u^+)^{p-1}hdx|\nonumber\\
&\leq&C'(\int_{\mathbb{R}^N}(|u_n(\cdot+y_n)|^p+|u|^p)dx)^{\frac{p-1}{p}}(\int_{\mathbb{R}^N}|K_*(x+y_n)-\mu|^p\cdot|h|^pdx)^{\frac{1}{p}}
\nonumber\\
&\leq&C(\int_{\mathbb{R}^N}|K_*(x+y_n)-\mu|^p\cdot|h|^pdx)^{\frac{1}{p}}\rightarrow
0\nonumber
\end{eqnarray}
where $C'$ and $C$ are positive constants independent of $n$ and
$h$. This together with (\ref{hftryfgggf}) and (\ref{mxxxcxcxx})
yields
\begin{eqnarray}\label{kkhjbnnb88yuh}
\langle J'(u_n), h(\cdot-y_n)\rangle=\langle J'_\theta(u),
h\rangle+o(1)
\end{eqnarray}
Then  by the assumption $J'(u_n)\rightarrow 0$ in
$H^{-1}(\mathbb{R}^N)$, we get $\langle J'_\theta(u), h\rangle=0$,
$\forall h\in H^1(\mathbb{R}^N).$ Therefore, $J'_\theta(u)=0.$

 3). From the definition of $v_n$,
\begin{eqnarray}\label{k77rtrgfg}
||v_n||^2_A=||u_n-u(\cdot-y_n)||^2_A=||u_n||^2_A+||u(\cdot-y_n)||^2_A-2(u_n,u(\cdot-y_n))_A.
\end{eqnarray}
By the definition of the norm $||\cdot||_A$ (see
(\ref{uutyghhgg})), we have
\begin{eqnarray}\label{i775656tygg}
||u(\cdot-y_n)||^2_A &=&\int_{\mathbb{R}^N}|\nabla
u(\cdot-y_n)|^2dx+(\frac{b^2}{4}-b)\int_{\mathbb{R}^N}\frac{|x\cdot\nabla
u(\cdot-y_n)|^2}{|x|^2}dx\nonumber\\
&&+\int_{\mathbb{R}^N}V_*(x)|u(\cdot-y_n)|^2dx\nonumber\\
&=&\int_{\mathbb{R}^N}|\nabla
u|^2dx+(\frac{b^2}{4}-b)\int_{\mathbb{R}^N}\frac{|(\frac{x}{|y_n|}+\frac{y_n}{|y_n|})\cdot\nabla
u|^2}{|\frac{x}{|y_n|}+\frac{y_n}{|y_n|}|^2}dx\nonumber\\
&&+\int_{\mathbb{R}^N}V_*(x+y_n)|u|^2dx.
\end{eqnarray}
Since $\nabla u\in L^2(\mathbb{R}^N)$ and
$\frac{x}{|y_n|}+\frac{y_n}{|y_n|}\rightarrow\theta$ a.e. on
$\mathbb{R}^N$, using the Lebesgue convergence theorem, we get
that
\begin{eqnarray}\label{mmvnvbv2q}
\int_{\mathbb{R}^N}\frac{|(\frac{x}{|y_n|}+\frac{y_n}{|y_n|})\cdot\nabla
u|^2}{|\frac{x}{|y_n|}+\frac{y_n}{|y_n|}|^2}dx\rightarrow\int_{
\mathbb{R}^N}|\theta\cdot\nabla u|^2dx,\ n\rightarrow\infty.
\end{eqnarray}
By (\ref{uuthgbbvb}), $(\ref{i775656tygg})$ and (\ref{mmvnvbv2q}),
we get that
\begin{eqnarray}\label{hgyyt6656rt}
||u(\cdot-y_n)||^2_A &=&\int_{\mathbb{R}^N}|\nabla
u|^2dx+(\frac{b^2}{4}-b)\int_{\mathbb{R}^N}|\theta\cdot\nabla
u|^2dx+a\int_{\mathbb{R}^N}|u|^2dx+o(1)\nonumber\\
&=&||u||^2_\theta+o(1).
\end{eqnarray}
Combining (\ref{k77rtrgfg}), (\ref{hgyyt6656rt}) and
(\ref{mxxxcxcxx}) leads to
\begin{eqnarray}\label{uut6657tytt}
||v_n||^2_A=||u_n||^2_A-||u||^2_\theta+o(1).
\end{eqnarray}

Note that
\begin{eqnarray}\label{iutygh77tyg}
&&\int_{\mathbb{R}^N}K_*(x)(v^+_n)^pdx\nonumber\\
&=&\int_{\mathbb{R}^N}K_*(x+y_n)((u_n(\cdot+y_n)-u)^+)^pdx\nonumber\\
&=&\int_{\mathbb{R}^N}((K^{\frac{1}{p}}_*(x+y_n)u_n(\cdot+y_n)-K^{\frac{1}{p}}_*(x+y_n)
u)^+)^pdx.
\end{eqnarray}
  We obtain from Lemma \ref{honhtton} that
\begin{eqnarray}\label{igythrgrtf66t}
&&\int_{\mathbb{R}^N}((K^{\frac{1}{p}}_*(x+y_n)u_n(\cdot+y_n)-K^{\frac{1}{p}}_*(x+y_n)
u)^+)^pdx\nonumber\\
&=&\int_{\mathbb{R}^N}((K^{\frac{1}{p}}_*(x+y_n)u^+_n(\cdot+y_n))^pdx
-\int_{\mathbb{R}^N}((K^{\frac{1}{p}}_*(x+y_n)u^+)^pdx+o(1)\nonumber\\
&=&\int_{\mathbb{R}^N}K_*(x)(u^+_n)^pdx-\int_{\mathbb{R}^N}K_*(x+y_n)(u^+)^pdx+o(1).
\end{eqnarray}
 By Lemma
\ref{nnvbf5rtf},
\begin{eqnarray}\label{kkg77t65ygf}
\int_{\mathbb{R}^N}K_*(x+y_n)(u^+)^pdx=\mu\int_{\mathbb{R}^N}(u^+)^pdx+o(1).
\end{eqnarray}
Combining $(\ref{iutygh77tyg})-(\ref{kkg77t65ygf})$ yields
\begin{eqnarray}\label{iiiiaaa}
\int_{\mathbb{R}^N}K_*(x)(v^+_n)^pdx=\int_{\mathbb{R}^N}K_*(x)(u^+_n)^pdx-\mu\int_{\mathbb{R}^N}(u^+)^pdx+o(1).
\end{eqnarray}
Combining (\ref{uut6657tytt}), (\ref{iiiiaaa}) and the assumption
$J(u_n)\rightarrow c$ leads to
\begin{eqnarray}\label{ookghbny7}
J(v_n)=J(u_n)-J_\theta(u)+o(1)=c-J_\theta(u)+o(1).\nonumber
\end{eqnarray}

4). For $h\in H^1(\mathbb{R}^N)$,
\begin{eqnarray}\label{kmnfgd55er}
\langle J'(v_n),
h\rangle&=&(v_n,h)_A-\int_{\mathbb{R}^N}K_*(x)(v_n^+)^{p-1}hdx\nonumber\\
&=&(u_n,
h)_A-(u(\cdot-y_n),h)_A-\int_{\mathbb{R}^N}K_*(x)(v_n^+)^{p-1}hdx.
\end{eqnarray}
We shall give the limits for $(u(\cdot-y_n),h)_A$ and
$\int_{\mathbb{R}^N}K_*(x)(v_n^+)^{p-1}hdx$ as
$n\rightarrow\infty.$

First, as (\ref{ii9945554es}), we have
\begin{eqnarray}\label{iugdvcxx1s}
&&(u(\cdot-y_n),h)_A\nonumber\\
&=&\int_{\mathbb{R}^N}\nabla u\nabla
h(\cdot+y_n)dx+a\int_{\mathbb{R}^N}
u\cdot h(\cdot+y_n)dx\nonumber\\
&&+\int_{\mathbb{R}^N}(V_*(x+y_n)-a)u\cdot
h(\cdot+y_n)dx\nonumber\\
&&+(\frac{b^2}{4}-b)\int_{\mathbb{R}^N}\Big(\frac{\frac{x}{|y_n|}+\frac{y_n}{|y_n|}}{|\frac{x}{|y_n|}+\frac{y_n}{|y_n|}|}\cdot\nabla
u\Big)\Big(\frac{\frac{x}{|y_n|}+\frac{y_n}{|y_n|}}{|\frac{x}{|y_n|}+\frac{y_n}{|y_n|}|}\cdot\nabla
h(\cdot+y_n)\Big)dx.\nonumber
\end{eqnarray}
By the H\"older inequality and (\ref{uuthgbbvb}), we get that if
$||h||\leq 1,$ then
\begin{eqnarray}\label{lmmcv66ftg955}
&&|\int_{\mathbb{R}^N}(V_*(x+y_n)-a)u\cdot
h(\cdot+y_n)dx|\nonumber\\
&\leq& (\int_{\mathbb{R}^N}|V_*(x+y_n)-a|\cdot u^2dx
)^{1/2}(\int_{\mathbb{R}^N}|V_*(x)-a|h^2dx )^{1/2}\nonumber\\
&\leq&C(\int_{\mathbb{R}^N}|V_*(x+y_n)-a|\cdot u^2dx
)^{1/2}\rightarrow 0,\ n\rightarrow\infty.\nonumber
\end{eqnarray}
 Thus, as
$n\rightarrow\infty,$
 \begin{eqnarray}\label{odgcbbcv9i}
\int_{\mathbb{R}^N}(V_*(x+y_n)-a)u\cdot
h(\cdot+y_n)dx=o(1)\nonumber
\end{eqnarray}
holds uniformly for $||h||\leq 1$. Moreover, a similar argument as
the proof of $(\ref{mmgyrttr6})$  yields that, as
$n\rightarrow\infty$,
\begin{eqnarray}
&&\int_{\mathbb{R}^N}\Big(\frac{\frac{x}{|y_n|}+\frac{y_n}{|y_n|}}{|\frac{x}{|y_n|}+\frac{y_n}{|y_n|}|}\cdot\nabla
u\Big)\Big(\frac{\frac{x}{|y_n|}+\frac{y_n}{|y_n|}}{|\frac{x}{|y_n|}+\frac{y_n}{|y_n|}|}\cdot\nabla
h(\cdot+y_n)\Big)dx\nonumber\\
&=&\int_{\mathbb{R}^N}(\theta\cdot\nabla u)(\theta\cdot\nabla
h(\cdot+y_n))dx+o(1)\nonumber
\end{eqnarray}
holds uniformly for $||h||\leq 1$. Therefore, as
$n\rightarrow\infty,$
\begin{eqnarray}\label{kk77tytyt}
(u(\cdot-y_n),h)_A=(u,h(\cdot+y_n))_\theta+o(1)
\end{eqnarray}
holds uniformly for $||h||\leq 1.$

Second, from  $u_n(\cdot+y_n)\rightharpoonup u$ in
$H^1(\mathbb{R}^N)$ and Lemma \ref{ooh8y77yuyy}, we deduce that,
as $n\rightarrow\infty,$
\begin{eqnarray}\label{mmkiswaas}
&&|\int_{\mathbb{R}^N}K_*(x+y_n)((u_n(\cdot+y_n)-u)^+)^{p-1}h(\cdot+y_n)dx\nonumber\\
&&\quad-\int_{\mathbb{R}^N}K_*(x+y_n)((u^+_n(\cdot+y_n))^{p-1}h(\cdot+y_n)dx\nonumber\\
&&\quad+\int_{\mathbb{R}^N}K_*(x+y_n)(u^+)^{p-1}h(\cdot+y_n)dx|\rightarrow
0
\end{eqnarray}
holds uniformly for $||h||\leq 1.$  By the H\"older inequality,
(\ref{nnvbf66r5}) and Lemma \ref{nnvbf5rtf}, we get that, if
$||h||\leq 1$, then
\begin{eqnarray}\label{jjnvyy5trtr}
&&|\int_{|x|>R}K_*(x+y_n)(u^+)^{p-1}h(\cdot+y_n)dx|\nonumber\\
&\leq&
(\int_{|x|>R}K_*(x+y_n)(u^+)^{p}dx)^{\frac{p-1}{p}}(\int_{|x|>R}K_*(x+y_n)|h|^{p}dx)^{1/p}\nonumber\\
&\leq&
C(\int_{|x|>R}K_*(x+y_n)(u^+)^{p}dx)^{\frac{p-1}{p}}\rightarrow
0,\ R\rightarrow\infty
\end{eqnarray}
 By Lemma
\ref{nnvbf5rtf}, we get that, for every $R>0,$  as
$n\rightarrow\infty,$
\begin{eqnarray}\label{jgjnvyy5trtr}
&&\sup_{||h||\leq 1}|\int_{|x|\leq R}(K_*(x+y_n)-\mu)(u^+)^{p-1}h(\cdot+y_n)dx|\nonumber\\
&\leq&
\sup_{||h||\leq 1}(\int_{|x|\leq R}|K_*(x+y_n)-\mu|(u^+)^{p}dx)^{\frac{p-1}{p}}(\int_{ \mathbb{R}^N}|K_*(x)-\mu|\cdot|h|^{p}dx)^{1/p}\nonumber\\
&\leq& C(\int_{|x|\leq
R}|K_*(x+y_n)-\mu|(u^+)^{p}dx)^{\frac{p-1}{p}}\rightarrow 0.
\end{eqnarray}
Combining (\ref{jjnvyy5trtr}) and (\ref{jgjnvyy5trtr}) yields that
\begin{eqnarray}\label{mmvnvbyyft6}
\int_{\mathbb{R}^N}K_*(x+y_n)(u^+)^{p-1}h(\cdot+y_n)dx-\mu\int_{\mathbb{R}^N}(u^+)^{p-1}h(\cdot+y_n)dx
\rightarrow 0
\end{eqnarray}
holds uniformly for $||h||\leq 1.$ Then by (\ref{mmkiswaas}),
(\ref{mmvnvbyyft6}) and \begin{eqnarray}
\int_{\mathbb{R}^N}K_*(x)(v_n^+)^{p-1}hdx=\int_{\mathbb{R}^N}K_*(x+y_n)((u_n(\cdot+y_n)-u)^+)^{p-1}hdx,\nonumber
\end{eqnarray} we get that, as $n\rightarrow\infty,$
\begin{eqnarray}\label{mmkiswafa}
&&|\int_{\mathbb{R}^N}K_*(x)(v_n^+)^{p-1}hdx-\int_{\mathbb{R}^N}K_*(x)(u^+_n)^{p-1}hdx
+\mu\int_{\mathbb{R}^N}(u^+)^{p-1}h(\cdot+y_n)dx|\nonumber\\
&&\rightarrow 0
\end{eqnarray}
holds uniformly for $||h||\leq 1.$

 Finally, combining (\ref{kmnfgd55er}),
(\ref{kk77tytyt}) and (\ref{mmkiswafa}) leads to
\begin{eqnarray}
\langle J'(v_n), h\rangle-\langle J'(u_n), h\rangle+\langle
J'_\theta(u), h(\cdot+y_n)\rangle\rightarrow 0\nonumber
\end{eqnarray}
holds uniformly for $||h||\leq 1.$ This together with the fact
that $J'_\theta(u)=0$ and $J'(u_n)\rightarrow 0$ in
$H^{-1}(\mathbb{R}^N)$ yields $J'(v_n)\rightarrow 0$ in
$H^{-1}(\mathbb{R}^N)$.\hfill$\Box$

\bigskip

\noindent{\bf Proof of Theorem \ref{uut6ghvbbo}.} We divide  the
proof into two steps.

1). For $n$ big enough, we have
\begin{eqnarray}\label{jjf665trgg}
c+1+||u_n||\geq J(u_n)-p^{-1}\langle J'(u_n),
u_n\rangle=(\frac{1}{2}-\frac{1}{p})||u_n||^2_A.
\end{eqnarray}
As mentioned in section \ref{tte54rrer}, the norm $||\cdot||_A$
 is
equivalent to  the norm $||\cdot||$. Therefore, there exists a
constant $C>0$ such that $||u||_A\geq C||u||,$ $\forall u\in
H^1(\mathbb{R}^N)$. Then by (\ref{jjf665trgg}), there exists a
constant $C'>0$ such that for $n$ big enough,
$$c+1+||u_n||\geq C'||u_n||^2$$
 It follows that $||u_n||$ is bounded.

 2). Assume that $u_n\rightharpoonup u_0$ in $H^1(\mathbb{R}^N)$
 and $u_n\rightarrow u_0$ a.e. on $\mathbb{R}^N$. By Lemma \ref{iig8uuhn},
 $J'(u_0)=0$ and $u^1_n=u_n-u_0$ is such that
 \begin{eqnarray}\label{99rytttytcc}
& &||u^1_{n}||^{2}_A=||u_{n}||^{2}_A-||u_0||^{2}_A+o(1),\nonumber \\
& &J(u^1_{n})\rightarrow c-J(u),\\
& &J'(u^1_{n})\rightarrow 0  \ \mbox{in} \ H^{-1}(\mathbb{R}^N).
\nonumber
\end{eqnarray}
Let us define
$$\delta:=\overline{\lim}_{n\rightarrow\infty}\sup_{y\in\mathbb{R}^N}\int_{|x-y|\leq 1}|u^1_n|^2dx.$$
If $\delta=0,$ Lemma \ref{llg66ghyhhv} implies that
$K^{1/p}_*u^1_n\rightarrow 0$ in $L^p(\mathbb{R}^N)$. Since
$J'(u^1_n)\rightarrow 0$ in $H^1(\mathbb{R}^N)$, it follows that
$$||u^1_n||^2_A=\langle J'(u^1_n), u_n^1\rangle+\int_{\mathbb{R}^N}K_*(x)((u^1_n)^+)^pdx\rightarrow 0$$
and the proof is complete. If $\delta>0,$ we may assume the
existence of $\{y^1_n\}\subset \mathbb{R}^N$ such that
$$\int_{|x-y^1_n|\leq 1}|u^1_n|^2dx>\delta/2.$$
Let us define $v^1_n:=u^1_n(\cdot+y^1_n)$. We may assume that
$v^1_n\rightharpoonup u_1$ in $H^1(\mathbb{R}^N)$ and
$v^1_n\rightarrow u_1$ a.e. on $\mathbb{R}^N$. Since
$$\int_{|x|\leq 1}|v^1_n|^2dx>\delta/2$$
it follows from the Rellich Theorem that $$\int_{|x|\leq
1}|u^1|^2dx\geq\delta/2$$ and $u_1\neq 0.$ But
$u^1_n\rightharpoonup 0$ in $H^1(\mathbb{R}^N)$, so that
$\{|y^1_n|\}$ is unbounded. We may assume that
$|y^1_n|\rightarrow\infty$. Finally, by (\ref{99rytttytcc}) and
Lemma \ref{dondgdonguu}, there exists $\theta_1\in\mathbb{R}^N$
with $|\theta_1|=1$ such that $J'_{\theta_1}(u_1)=0$ and
$u^2_n:=u^1_n-u_1(\cdot-y^1_n)$ satisfies
\begin{eqnarray}\label{iiiiawwqwq}
& &||u^2_{n}||^{2}=||u^1_{n}||^{2}-||u_1||^{2}+o(1),\nonumber \\
& &J(u^2_{n})\rightarrow c-J_{\theta_1}(u_1), \nonumber\\
& &J'(u^2_{n})\rightarrow 0  \ \mbox{in} \ H^{-1}(\mathbb{R}^N).
\nonumber
\end{eqnarray}
Moreover,  Lemma \ref{iighjgytyf} implies that
 $$J_{\theta_1}(u_1)\geq(\frac{1}{2}-\frac{1}{p})\mu^{-\frac{2}{p-2}}S^{\frac{p}{p-2}}.$$

Iterating the above procedure we construct sequences
$\{\theta_l\}$, $\{u_l\}$ and $\{y^l_n\}$. Since for every $l,$
$J_{\theta_l}(u_l)\geq(\frac{1}{2}-\frac{1}{p})\mu^{-\frac{2}{p-2}}S^{\frac{p}{p-2}},$
the iteration must terminate at some finite index $k$. This
finishes the proof of this theorem.\hfill$\Box$

\bigskip

The following corollary is a direct consequence of Theorem
\ref{uut6ghvbbo} and Lemma \ref{iighjgytyf}. It implies that the
functional $J$ satisfies $(PS)_c$ condition if
$c<(\frac{1}{2}-\frac{1}{p})\mu^{-\frac{2}{p-2}}S^{\frac{p}{p-2}}.$

\begin{corollary}\label{kkg77t6yy}
Under the assumptions  $\bf (A_1)$ and $\bf (A_2)$,  any sequence
$\{u_n\}\subset H^1(\mathbb{R}^N)$ such that
$$ J(u_n)\rightarrow c<(\frac{1}{2}-\frac{1}{p})\mu^{-\frac{2}{p-2}}S^{\frac{p}{p-2}},\ J'(u_n)\rightarrow0
\ \mbox{in}\ H^{-1}(\mathbb{R}^N)$$ contains a convergent
subsequence.
\end{corollary}

\section{Proof of Theorem \ref{ll7trggfgf}. }
Recall that the critical points of $J$ are nonnegative solutions
of (\ref{o8867yyttggd}).  By Corollary \ref{kkgn775trff}, to prove
Eq.(\ref{h777arrrrrr}) has a positive solution, it suffices to
prove that  $J$ has a nontrivial critical point. And by Corollary
\ref{kkg77t6yy}, it suffices to apply the classical mountain pass
theorem (see,  e.g.,  \cite[Theorem 1.15]{Willem}) to $J$ with the
mountain pass value
$c<(\frac{1}{2}-\frac{1}{p})\mu^{-\frac{2}{p-2}}S^{\frac{p}{p-2}}.$

By the assumption (\ref{ofbbv66ft}) and Lemma \ref{fddretrdfdf},
there exists a nonnegative $u_0\in
H^1(\mathbb{R}^N)\setminus\{0\}$ such that
\begin{eqnarray}\label{wsofbbv66ft}
&&\frac{||u_0||_A^2}{(\int_{\mathbb{R}^N}K_*(x)u_0^pdx)^{2/p}}<(1-b/2)^{\frac{p-2}{p}}\mu^{-\frac{2}{p}}S_p=
\mu^{-\frac{2}{p}}S.\nonumber
\end{eqnarray}
We obtain
\begin{eqnarray}\label{yrttftt63}
0<\max_{t\geq 0}J(tu_0)&=&\max_{t\geq
0}\Big(\frac{t^2}{2}||u_0||^2_A-\frac{t^p}{p}\int_{\mathbb{R}^N}K_*(x)(u^+_0)^pdx\Big)\nonumber\\
&=&(\frac{1}{2}-\frac{1}{p})\Big(||u_0||^2_A\Big/(\int_{\mathbb{R}^N}K_*(x)u_0^pdx)^{2/p}\Big)^{\frac{p}{p-2}}\nonumber\\
&<&(\frac{1}{2}-\frac{1}{p})\mu^{-\frac{2}{p-2}}S^{\frac{p}{p-2}}.
\end{eqnarray}
By (\ref{nnvbf66r5}),
\begin{eqnarray}
J(u)\geq \frac{1}{2}||u||_A^2-\frac{C^p}{p}||u||^p_A.\nonumber
\end{eqnarray}
 Therefore, there exists $r>0$
such that $$b:=\inf_{||u||_A=r}J(u)>0=J(0).$$ Moreover, there
exists $t_0>0$ such that $||t_0u_0||_A>r$ and $J(t_0u_0)<0.$ It
follows from (\ref{yrttftt63})  that $$\max_{t\in [
0,1]}J(tt_0u_0)<(\frac{1}{2}-\frac{1}{p})\mu^{-\frac{2}{p-2}}S^{\frac{p}{p-2}}.$$
By   Corollary \ref{kkg77t6yy} and the mountain pass theorem (see
\cite[Theorem 1.15]{Willem}), $J$ has a critical value $c$ such
that $b\leq
c<(\frac{1}{2}-\frac{1}{p})\mu^{-\frac{2}{p-2}}S^{\frac{p}{p-2}}$
and Eq.(\ref{o8867yyttggd}) has a positive solution $v\in
H^1(\mathbb{R}^N)$. Then by Theorem \ref{kkgn775trff}, the
function $u$ defined by (\ref{lululuggftrr}) is a positive
solution of (\ref{h777arrrrrr}). To  complete the proof, it
suffices to prove that $u\in E$. Using the divergence theorem,
Lemma \ref{jjfgr665t2w} and (\ref{mmmmmazx}), we get that
\begin{eqnarray}\label{ooodgcbbcv}
&&\int_{\mathbb{R}^N}|\nabla
u|^2dx\nonumber\\
&=&-\int_{\mathbb{R}^N}u\triangle
u dx\nonumber\\
&=&-\int_{\mathbb{R}^N}u\cdot|y|^{-\frac{b(N+2)}{2(2-b)}}\Big(\sum_{i,j=1}^N\frac{\partial}{\partial
y_j}\Big(A_{ij}(y)\frac{\partial v}{\partial
y_i}\Big)-\frac{C_{b}}{|y|^2}v\Big)dx\nonumber\\
&=&-\int_{\mathbb{R}^N}|x|^{-\frac{b}{4}(N-2)}v(|x|^{-\frac{b_*}{2}}x)\cdot|y|^{-\frac{b(N+2)}{2(2-b)}}\Big(\sum_{i,j=1}^N\frac{\partial}{\partial
y_j}\Big(A_{ij}(y)\frac{\partial v}{\partial
y_i}\Big)-\frac{C_{b}}{|y|^2}v\Big)dx\nonumber\\
&=&-\int_{\mathbb{R}^N}
v\cdot\Big(\sum_{i,j=1}^N\frac{\partial}{\partial
y_j}\Big(A_{ij}(y)\frac{\partial v}{\partial
y_i}\Big)-\frac{C_{b}}{|y|^2}v\Big)dy\nonumber\\
&=&\int_{\mathbb{R}^N}\Big(|\nabla v|^2+\frac{|x\cdot\nabla
v|^2}{|x|^2}+\frac{C_b}{|x|^2}v^2\Big)dy.\nonumber
\end{eqnarray}
Moreover, by Lemma \ref{jjfgr665t2w} and (\ref{mmmmmazx}), we get
that
\begin{eqnarray}
\int_{\mathbb{R}^N}V(x)u^2dx=\int_{\mathbb{R}^N}
V(x)|x|^{-\frac{b}{2}(N-2)}v^2(|x|^{-\frac{b}{2}}x)dx=\int_{\mathbb{R}^N}
V_*(y)v^2dy.
\end{eqnarray}
Therefore, $||u||^2_E=||v||^2_A<\infty.$ \hfill$\Box$

\bigskip

%{\noindent \bf Acknowledgements.}  Shaowei Chen was supported by
%Science Foundation of Huaqiao University.

\end{document}